\documentclass[a4paper,12pt,draft]{article}
\usepackage{amsmath,amssymb,amsthm,amsfonts}

\newcommand{\expr}[1]{\left( #1 \right)}

\newcommand{\A}{\mathcal{A}}
\newcommand{\R}{\mathbf{R}}
\newcommand{\Z}{\mathbf{Z}}
\newcommand{\Rd}{{\mathbf{R}^d}}
\newcommand{\Rt}{{\mathbf{R}^2}}
\newcommand{\ind}{\mathbf{1}}

\newcommand{\interior}{int}

\newtheorem{thm}{Theorem}[section]

\newtheorem{lem}[thm]{Lemma}
\newtheorem{prop}[thm]{Proposition}
\newtheorem{cor}[thm]{Corollary}
\newtheorem{remark}[thm]{Remark}

\theoremstyle{definition}

\numberwithin{equation}{section}

\DeclareMathOperator{\dist}{dist}
\DeclareMathOperator{\diam}{diam}

\title{Spectral gap for stable process on convex planar double symmetric domains}

\author{Bart{\l}omiej Dyda and Tadeusz Kulczycki\\
\small{Institute of Mathematics and Computer Science}\\
\small{Wroc{\l}aw University of Technology}\\
\small{Wybrze\.ze Wyspia\'nskiego 27, 50-370 Wroc{\l}aw, Poland} }

\begin{document}
\sloppy \footnotetext{
%tk w ogole to wywalilem  
%\emph{2000 Mathematics Subject
%Classification:} Primary 31B25, 60J45 ???.
The first named author was supported by KBN grant 1 P03A 026 29
and RTN Harmonic Analysis and  Related Problems, contract 
HPRN-CT-2001-00273-HARP\\
The second named author was supported by KBN 
grant 1 P03A 020 28
and RTN Harmonic Analysis and  Related Problems, contract 
HPRN-CT-2001-00273-HARP.}

\maketitle

\begin{abstract}
We study the semigroup of the symmetric
$\alpha$-stable process in bounded domains in $\Rt$. We obtain
a variational formula for
the spectral gap, i.e. the difference between two first eigenvalues
of the generator of this semigroup. This variational formula allows us to obtain lower bound estimates of the
 spectral gap for convex planar domains which are symmetric with respect to both coordinate axes. For rectangles, using "midconcavity" of  the first eigenfunction \cite{BKM}, we obtain sharp upper and lower bound estimates of the spectral gap.
\end{abstract}

\section{Introduction} \label{sec:introduction}
 
In recent years many results have been obtained in spectral theory of semigroups of symmetric $\alpha$-stable processes $\alpha \in (0,2)$ in bounded domains in $\Rd$, see \cite{BLM}, \cite{M}, \cite{BK1}, \cite{Db}, \cite{DbM}, \cite{CS3}, \cite{CS4}, \cite{BKM}. One of the most interesting problems in spectral theory of such semigroups is a spectral gap estimate i.e. the estimate of $\lambda_2 - \lambda_1$ the difference between two first eigenvalues of the generator of this semigroup. Such estimate is a natural generalisation of the same problem for the semigroup of Brownian motion killed on exiting a bounded domain, which generator is Dirichlet Laplacian. In this classical case, for Brownian motion, spectral gap estimates have been widely studied see e.g \cite{SWYY}, \cite{YZ}, \cite{L}, \cite{Sm}, \cite{Da}, \cite{BM}. When a bounded domain is convex there have been obtained sharp lower-bound estimates of the spectral gap.
 
In the case of the semigroup of symmetric $\alpha$-stable processes $\alpha \in (0,2)$ very little is known about the spectral gap estimates. In one dimensional case when a domain is just an interval spectral gap estimates follow from results from \cite{BK1} ($\alpha = 1$) and \cite{CS3} ($\alpha > 1$). The only results for dimension greater than one have been obtained for the Cauchy process i.e. $\alpha = 1$ \cite{BK2}, \cite{BK3}. Such results have been obtained using the deep connection between the eigenvalue problem for the Cauchy process and a boundary value problem for the Laplacian in one dimension higher, known as the mixed Steklov problem.

The aim of this paper is to generalise spectral gap estimates obtained for the Cauchy process ($\alpha = 1$) for all $\alpha \in (0,2)$. Before we describe our results in more detail let us recall definitions and basic facts.
 
Let $X_t$ be a symmetric $\alpha$-stable process in $\Rd$, $\alpha \in
(0,2]$. This is a process with independent and stationary increments
and characteristic function 
$E^0 e^{i \xi X_t} = e^{-t |\xi|^{\alpha}}$, $\xi \in \Rd$, $ t > 0$. We will use 
$E^x$, $P^x$ to denote the expectation and probability of this process
starting at $x$, respectively.  By $p(t,x,y) = p_{t}(x-y)$ we
will denote the transition density of this process. That is, 
$$P^{x}(X_t \in
B)= \int_{B} p(t,x,y)\, dy.
$$
When $\alpha = 2$ the process $X_t$ is just the Brownian motion in $\Rd$ 
running at twice
the speed. That is, if $\alpha = 2$ then
\begin{equation}
\label{gaussian}
p(t,x,y) = \frac{1}{(4 \pi t)^{d/2}} e^{\frac{-|x-y|^2}{4t}},\quad t > 0, \, \, x,y \in \Rd.
\end{equation}

It is well known that for $\alpha \in (0,2)$ we have $p_t(x) = t^{-d/\alpha} p_1(t^{-1/\alpha} x)$, $t > 0$, $x \in \Rd$ and 
$$
p_t(x) = t^{-d/\alpha} p_1(t^{-1/\alpha} x) \le t^{-d/\alpha} p_1(0) = t^{-d/\alpha} M_{d,\alpha}, \quad t > 0, \,  x \in \Rd,
$$
where
\begin{equation}
\label{Mdalpha}
M_{d,\alpha} = \frac{1}{(2 \pi)^{d}} \int_{\Rd} e^{- |x|^{\alpha}} \, dx.
\end{equation}

%tk tu usunalem ???
It is also well known that
\begin{equation}
\label{limptxy}
\lim_{t \to 0^+} \frac{p(t,x,y)}{t} =
\frac{\A_{d,-\alpha}}{|x - y|^{d + \alpha}},
\end{equation}
where
\begin{equation}
\label{Adalpha}
\A_{d,\gamma} = \Gamma((d - \gamma)/2)/(2^{\gamma} \pi^{d/2} |\Gamma(\gamma/2)|).
\end{equation}

Our main concern  in this paper are the eigenvalues of the semigroup of the process
$X_t$ killed upon  leaving a domain. Let $D \subset \Rd$ be a bounded connected domain
and $\tau_{D} =
\inf\{t
\ge 0: X_t
\notin D\}$ be the first exit time of $D$.   By $\{P_t^{D}\}_{t \ge 0}$ we denote the semigroup
on
$L^2(D)$ of $X_t$ killed upon exiting $D$.  That is, 
 
$$
P_{t}^{D}f(x) = E^{x}(f(X_{t}), \tau_{D} > t), \quad x \in D, \, \, t > 0, \, \, f \in L^2(D).
$$
The semigroup has transition densities $p_D(t,x,y)$ satisfying 
$$
P_{t}^{D}f(x) = \int_{D} p_{D}(t,x,y) f(y) \, dy.
$$
The  kernel $p_{D}(t,x,y)$  is strictly positive
symmetric and
$$
p_{D}(t,x,y) \le p(t,x,y) \le M_{d,\alpha} \, t^{-d/\alpha}, \quad x,y \in D, \, \, t>0.
$$
The fact that $D$ is bounded implies that for any $t > 0$ the operator $P_t^D$ maps $L^2(D)$
into $L^{\infty}(D)$. From the
general theory of
 semigroups (see 
\cite{Da1})  it follows that there exists 
 an   orthonormal basis of eigenfunctions
$\{\varphi_n\}_{n =1}^{\infty}$ for
$L^2(D)$ and corresponding eigenvalues
$\{\lambda_n\}_{n = 1}^{\infty}$ satisfying
$$0<\lambda_1<\lambda_2\leq \lambda_3\leq \dots$$
  with $\lambda_n\to\infty$
as
$n\to\infty$. That is,  the pair $\{\varphi_n, \lambda_n\}$ satisfies
\begin{equation}
\label{cauchyproblem}
P_{t}^{D}\varphi_{n}(x) = e^{-\lambda_{n} t} \varphi_{n}(x), \quad x \in D, \,\,\, t > 0.
\end{equation}
The eigenfunctions $\varphi_n$ are continuous and bounded on $D$. In addition, 
$\lambda_1$ is simple and the corresponding eigenfunction
$\varphi_1$, often called the ground state eigenfunction, is strictly
positive on
$D$.  
By scaling we have for $\beta>0$
\begin{equation}\label{gapscaling}
 \lambda_n(\beta D)=\beta^{-\alpha} \lambda_n(D).
\end{equation}
 For more general properties of
  the semigroups $\{P_t^D\}_{t \ge 0}$, see
\cite{G1},
\cite{BG2}, \cite{CS1}.

It is well known (see \cite{B}, \cite{CS1}, \cite{CS2}, \cite{K})
that if $D$ is a bounded connected Lipschitz domain and $\alpha = 2$, or that if 
$D$ is a bounded connected domain for $0 < \alpha < 2$, then
$\{P_t^D\}_{t \ge 0}$ is intrinsically ultracontractive.  Intrinsic ultracontractivity is a 
remarkable property with many consequences.  It implies, in particular, that  
$$
\lim_{t \to \infty} \frac{e^{\lambda_1 t} p_D(t,x,y)}{\varphi_1(x) \varphi_1(y)} = 1, 
$$
uniformly in both variables $x, y\in D$.  In addition, the rate of convergence is given by 
the spectral gap $\lambda_2-\lambda_1$. That is, 
for any $t \ge 1$ we have
\begin{equation}
\label{spectralgap}
e^{-(\lambda_2 - \lambda_1) t} \le 
\sup_{x,y \in D} \left|\frac{e^{\lambda_1 t}p_D(t,x,y)}{ \varphi_1(x) \varphi_1(y)} - 1
 \right| \le C(D,\alpha) e^{-(\lambda_2 - \lambda_1)t}.
\end{equation}
The proof of this  for $\alpha=2$ may be found in \cite{Sm}.  
The proof in our setting is exactly the same.  

Our first step in studying the spectral gap for $\alpha \in (0,2)$ is the following variational characterisation of $\lambda_2 - \lambda_1$.

By $L^2(D,\varphi_1^2)$ we denote the $L^2$ space of functions with
the inner product $(f,g)_{L^2(D,\varphi_1^2)} = \int_D f(x) g(x)
\varphi_1^2(x) \, dx$.

\begin{thm}
\label{variational}
We have
\begin{equation}
\label{variational1}
\lambda_2 - \lambda_1
= \inf_{f \in \mathcal{F}} \frac{\A_{d,-\alpha}}{2} \int_D \int_D
\frac{(f(x) - f(y))^2}{|x - y|^{d + \alpha}} \varphi_1(x) \varphi_1(y) \, dx \, dy,
\end{equation}
where
$$
\mathcal{F} = \{f \in L^2(D,\varphi_1^2): \int_D f^2(x) \varphi_1^2(x) \, dx = 1, \, \, \, \int_D f(x) \varphi_1^2(x) \, dx = 0\}
$$
and $\A_{d,-\alpha}$ is given by (\ref{Adalpha}). Moreover the infimum is achieved for $f = \varphi_2/\varphi_1$.
\end{thm}
The idea of the proof is based on considering a new semigroup $\{T_t^D\}_{t \ge 0}$ of the stable process conditioned to remain forever in $D$. The proof of Theorem~\ref{variational} is in Section 2.

In the classical case, for Brownian motion, when a dimension is greater than one, the simplest domain where the spectral gap can be explicitly calculated is a rectangle. Let us recall that in this classical case $\{\varphi_n\}_{n = 1}^{\infty}$, $\{\lambda_n\}_{n = 1}^{\infty}$ are of course eigenfunctions and eigenvalues of Dirichlet Laplacian. Therefore, when (say) $D = (-a,a) \times (-b,b)$, $a \ge b > 0$ then
$$
\varphi_1(x_1,x_2) = (1/\sqrt{2 a b}) \cos(\pi x_1/(2 a))
\cos(\pi x_2/(2 b)), 
$$
$$
\varphi_2(x_1,x_2) = (1/\sqrt{2 a b}) \sin(2 \pi x_1/(2 a))
\cos(\pi x_2/(2 b)),
$$
$\lambda_1 = \pi^2/(4 a^2) + \pi^2/(4 b^2)$, $\lambda_2 = 4\pi^2/(4 a^2) + \pi^2/(4 b^2)$ and hence
$\lambda_2 - \lambda_1 = 3 \pi^2/(4 a^2)$.

Although the $\alpha$-stable process is generated by $-(-\Delta)^{\alpha/2}$, the generator of the killed $\alpha$-stable process on $D$ is however not equal to $-(-\Delta_D)^{\alpha/2}$ for the Dirichlet Laplacian $\Delta_D$ on $D$. So, both $\varphi_n$ and $\lambda_n$ are not explicit even for an interval or a rectangle. However, when $D$ is a rectangle, due to simple geometric properties of this set it is shown (\cite{BKM} Theorem~1.1) that the first eigenfunction $\varphi_1$ for any $\alpha \in (0,2]$ is "midconcave" and unimodal according to the lines parallel to the sides. This property and Theorem~\ref{variational} enables us to obtain sharp upper and lower bound estimates of the spectral gap for all $\alpha \in (0,2)$. The most complicated are lower bound estimates for $\alpha \in (1,2)$ and $\alpha = 1$. The main idea of the proof in these cases is contained in Lemmas~\ref{partition0} and \ref{partition}.

Below we present estimates of $\lambda_2 - \lambda_1$ for rectangles. The proof of this theorem is in Section 4. Let us point out that these estimates are sharp i.e. the upper and lower bound estimates have the same dependence on the length of the sides of the rectangle. Nevertheless, the numerical constants which appear in this theorem are far from being optimal.
\begin{thm}\label{estrectangle}
Let $D=(-a,a)\times(-b,b)$, where $a\geq b$. Then
\begin{itemize}
\item[(a)] We have
$$
  2\mathcal{A}_{2,-\alpha}^{-1} (
  \lambda_2-\lambda_1) \le 10^6 \cdot
  \left\{ \begin{array}{ll}
       \displaystyle
       \frac{2}{1-\alpha} \frac{b}{a^{1+\alpha}} & 
            \textrm{for $\alpha<1$,} \vspace{2mm}\\
    \displaystyle
      2 \log\left(1 + \frac{a}{b}\right) \frac{b}{a^2} & \textrm{for $\alpha=1$,}  \vspace{2mm}\\
\displaystyle
      \left( \frac{1}{2-\alpha}+\frac{1}{\alpha-1}\right) \frac{b^{2 - \alpha}}{a^2} & 
            \textrm{for $\alpha>1$.} 
\end{array} \right.
$$
\item[(b)]
We have
\begin{eqnarray*}
   2\mathcal{A}_{2,-\alpha}^{-1} (
  \lambda_2-\lambda_1) \geq 
\left\{ \begin{array}{ll}
  \displaystyle  
  \frac{b}{36\cdot 2^{1 + 2 \alpha}  a^{1+\alpha} }
  &\textrm{for $\alpha<1$,} \vspace{2mm}
  \\
  \displaystyle  
      10^{-9} \, \log\left(1 + \frac{a}{b}\right) \frac{b}{a^2}& \textrm{for $\alpha=1$,}  \vspace{2mm}\\
  \displaystyle 
  \frac{1}{33 \cdot 13^{1 + \alpha/2} \cdot 10^4}\, \frac{b^{2 - \alpha}}{a^2}
  & \textrm{for $\alpha>1$.} 
    \end{array}
  \right.
\end{eqnarray*}
\end{itemize}
\end{thm}

Let us note that for $\alpha = 1$ the following estimates have already been known $\lambda_2 - \lambda_1 \ge C b/a^2$, where $C = 10^{-7}$ (Corollary 1.1, \cite{BK3}). However, estimates from Theorem \ref{estrectangle} are more precise because we get an extra term $\log(a/b + 1)$, which gives a sharp dependence on the length of the sides of a rectangle.

\begin{remark}
\label{remark1}
The inequality 
$$
2\mathcal{A}_{2,-\alpha}^{-1} (
  \lambda_2-\lambda_1) \geq 
  \frac{b}{36\cdot 2^{\alpha}  (a + b)^{1+\alpha} }
$$
holds for all $\alpha \in (0,2)$. 

We have 
 $2 \mathcal{A}_{2,-\alpha}^{-1}
 = \alpha^{-2} 2^{3 - \alpha} \pi \Gamma^{-1}(\alpha/2) \Gamma(1 - \alpha/2)$.
 In particular we get for example 
$\lambda_2 - \lambda_1 \ge \frac{8 b}{10^{4}(a + b)^{3/2}}$ 
 for $\alpha = 1/2$,
  $\lambda_2 - \lambda_1 \ge \frac{b}{10^3 (a + b)^2}$
   for $\alpha = 1$, $\lambda_2 - \lambda_1 \ge \frac{8 b}{10^4 (a + b)^{5/2}}$ for $\alpha = 3/2$.
\end{remark}

Our next aim are lower bound estimates of the spectral gap for convex planar domains which are symmetric with respect to both coordinate axes.

In the classical case, for the Brownian motion, there are known sharp estimates for all bounded convex domains $D \subset \Rd$. We have $\lambda_2 - \lambda_1 > \pi^2/d_D^2$ where $d_D$ is the diameter of $D$ see e.g. \cite{L}, \cite{Sm}. Such results are obtained using the fact that the first eigenfunction is log-concave. For convex planar domains which are symmetric with respect to both coordinate axes even better estimates $\lambda_2 - \lambda_1 > 3 \pi^2/d_D^2$ are known, see \cite{Da}, \cite{BM} (such estimates are optimal, the lower bound is approached by this rectangles). These results follow from ratio inequalities for heat kernels.

Unfortunately in the case of symmetric $\alpha$-stable processes, $\alpha \in (0,2)$, we do not know whether the first eigenfunction is log-concave. Instead we use some of the ideas from \cite{BK3} where spectral gap estimates for the Cauchy process i.e. $\alpha = 1$ were obtained. Namely, we use the fact that the first eigenfunction is unimodal according to the lines parallel to coordinate axes and that it satisfies the appropriate Harnack inequality. Then we use similar techniques as for rectangles. As before in this proof the crucial role have Lemmas~\ref{partition0} and \ref{partition}. 

The properties of the first eigenfunction are obtained in Section 3 and the proof of lower bound estimates for the spectral gap is in Section 5. 
These estimates we present below in Theorem~\ref{gapconvexset}.
Let us point out that these estimates are sharp only for $\alpha > 1$, where we know that they cannot be improved because of the results for rectangles.

\begin{thm}\label{gapconvexset}
Let $D \subset \Rt$ be a bounded convex domain which is symmetric relative
to both coordinate axes. Assume that $[-a,a] \times [-b,b]$, $a \ge b$ is the smallest rectangle (with sides parallel to the coordinate axes) containing $D$. 
Then we have
$$
 2\mathcal{A}_{2,-\alpha}^{-1} (\lambda_2-\lambda_1)
  \geq  \frac{C \, b^{2 - \alpha}}{a^2},
$$
where
\begin{equation}
\label{constantC}
C = C(\alpha) = 10^{-9} 3^{\alpha - 4} 2^{- 2 \alpha -1} \left(4 + \frac{12 \Gamma(2/\alpha)}{\alpha (2 - \alpha) (1 - 2^{-\alpha})^{2/\alpha}}\right)^{- 2}.
\end{equation}
\end{thm}

Let us note that for $\alpha = 1$ such estimate has already been known with a better constant. In fact, Corollary 1.1 \cite{BK3} gives $\lambda_2 - \lambda_1 \ge C b/a^2$, where $C = 10^{-7}$.

There are still many open problems concerning the spectral gap for semigroups of symmetric stable processes $\alpha \in (0,2)$ in bounded domains $D \subset \Rd$. Perhaps the most interesting is the following. What is the best possible lower bound estimate for the spectral gap for arbitrary bounded convex domain $D \subset \Rd$? With this problem there are connected questions about the shape of the first eigenfunction $\varphi_1$. For example, is $\varphi_1$ log-concave or at least unimodal when $D$ is a convex bounded domain? There is also an unsolved problem concerning domains from
Theorem~\ref{gapconvexset}. Can one obtain for $\alpha \le 1$ lower bounds similar to these obtained for rectangles i.e. $\lambda_2 - \lambda_1 \ge C_{\alpha} \, b /a^{1 + \alpha}$ for $\alpha < 1$ and $\lambda_2 - \lambda_1 \ge C \, b \log(1 + a/b)/a^{2}$ for $\alpha = 1$?

\section{Variational formula} \label{sec:variationalformula}

In this section we prove Theorem~\ref{variational} -- the variational formula
for the spectral gap.

At first we need the following simple properties of the kernel $p_D(t,x,y)$.

\begin{lem}
\label{transition}
There exists a constant $c = c(d,\alpha)$ such that for any $t > 0$, $
x,y \in D$ we have 

\begin{equation}
\label{transition1}
p_D(t,x,y) \le p(t,x,y) \le \frac{c t}{|x - y|^{d + \alpha}}.
\end{equation}
For any $x,y \in D$, $x \ne y$ we have
\begin{equation}
\label{transition2}
\lim_{t \to 0^+} \frac{p_D(t,x,y)}{t} = 
\lim_{t \to 0^+} \frac{p(t,x,y)}{t} =
\frac{\A_{d,-\alpha}}{|x - y|^{d + \alpha}}.
\end{equation}
\end{lem}

\begin{proof}
These properties of $p_D(t,x,y)$ are rather well known. We recall some of the standard arguments.

The estimate $p(t,x,y) \le c t |x - y|^{-d-\alpha}$ follows e.g. from the scaling property $p(t,x,y) = t^{-d/\alpha} p_1((x-y)t^{-1/\alpha})$ and the inequality $p_1(z) \le c |z|^{-d-\alpha}$ \cite{Z}.
The equality on the right-hand side of (\ref{transition2}) is well known (see (\ref{limptxy})).

We know that $p_D(t,x,y) = p(t,x,y) - r_D(t,x,y)$ where
$$
r_D(t,x,y) = E^x(\tau_D < t; p(t - \tau_D, X(\tau_D),y)).
$$
By (\ref{transition1}) we get for $x,y \in D$, $t > 0$
\begin{eqnarray*}
\frac{1}{t} r_D(t,x,y) &=&
\frac{1}{t} E^x(\tau_D < t; p(t - \tau_D, X(\tau_D),y)) \\
&\le& \frac{1}{t} E^x\left(\tau_D < t; \frac{c t}{|y - X(\tau_D)|^{d +
      \alpha}}\right) \\
&\le& \frac{c P^{x}(\tau_D < t)}{(\delta_D(y))^{d + \alpha}},
\end{eqnarray*}
where $\delta_D(y) = \inf \{|z - y|: z \in \partial{D}\}$. It follows that $t^{-1} r_D(t,x,y) \to 0$ when $t \to 0^+$.
\end{proof}

Let 
$$
\tilde{p}_D(t,x,y) = 
\frac{e^{\lambda_1 t} p_D(t,x,y)}{\varphi_1(x) \varphi_1(y)}, \quad x,y \in D, \quad t > 0
$$
and 
$$
T_t^D f (x) = \int_D \tilde{p}_D(t,x,y) f(y) \varphi_1^2(y) \, dy,
\quad f \in L^2(D,\varphi_1^2), \quad t > 0.
$$
$\{T_t^D\}_{t \ge 0}$ is a semigroup in $L^2(D,\varphi_1^2)$. This is
the semigroup for the stable process conditioned to remain forever in
$D$ (see \cite{Sm} where the same semigroup is defined for Brownian
motion).

Let 
$$ 
\mathcal{E}(f,f) 
= \lim_{t \to 0^+} \frac{1}{t} (f - T_t^Df, f)_{L^2(D,\varphi_1^2)},
$$
for $f \in L^2(D,\varphi_1^2)$.

\begin{lem}
\label{Dirichlet}
For any $f \in L^2(D,\varphi_1^2)$ $\mathcal{E}(f,f)$ is well defined and we have
\begin{equation}
\label{Dirichlet0}
\mathcal{E}(f,f)
= \frac{\A_{d,-\alpha}}{2} \int_D \int_D
\frac{(f(x) - f(y))^2}{|x - y|^{d + \alpha}} \varphi_1(x) \varphi_1(y) \, dx \, dy.
\end{equation}
\end{lem}

\begin{proof}
\begin{eqnarray}
\nonumber
&& \mathcal{E}(f,f) = 
\lim_{t \to 0^+} \frac{1}{t} (f - T_t^Df, f)_{L^2(D,\varphi_1^2)} \\
\nonumber
&=& \lim_{t \to 0^+} \frac{1}{t} 
\int_D \left(f(x) - 
\int_D \frac{e^{\lambda_1 t} p_D(t,x,y)}{\varphi_1(x) \varphi_1(y)} f(y) \varphi_1^2(y) \, dy \right)
f(x) \varphi_1^2(x) \, dx \\
\label{Dirichlet1}
&=& \lim_{t \to 0^+} \frac{1}{t} 
\int_D \left(f(x) \varphi_1(x) - 
e^{\lambda_1 t} \int_D p_D(t,x,y) f(y) \varphi_1(y) \, dy \right)\\
\nonumber
&& \times
f(x) \varphi_1(x) \, dx. 
\end{eqnarray}
Note that 
$$
f(x) \varphi_1(x) = f(x) e^{\lambda_1 t} P_t^D \varphi_1(x) =
e^{\lambda_1 t} \int_D p_D(t,x,y) f(x) \varphi_1(y) \, dy.
$$
Hence (\ref{Dirichlet1}) is equal to 
\begin{eqnarray}
\nonumber
&& \lim_{t \to 0^+} \frac{1}{t} 
\int_D  e^{\lambda_1 t} \int_D p_D(t,x,y) (f(x) \varphi_1(y) - f(y) \varphi_1(y)) \, dy 
f(x) \varphi_1(x) \, dx \\
\label{Dirichlet2}
&=& \lim_{t \to 0^+} e^{\lambda_1 t} 
\int_D \int_D \frac{p_D(t,x,y)}{t} (f^2(x) - f(x) f(y)) \varphi_1(x) \varphi_1(y) \, dy \, dx.
\end{eqnarray}
Note that we can interchange the role of $x$ and $y$ in
(\ref{Dirichlet2}). Therefore by standard arguments (\ref{Dirichlet2}) is equal to 
\begin{equation}
\label{Dirichlet3}
\lim_{t \to 0^+} \frac{e^{\lambda_1 t}}{2} 
\int_D \int_D \frac{p_D(t,x,y)}{t} (f(x) - f(y))^2 \varphi_1(x) \varphi_1(y) \, dx \, dy.
\end{equation}
In view of $(\ref{transition2})$ in order to prove (\ref{Dirichlet0}) we need only to justify the interchange of the limit and the integral in (\ref{Dirichlet3}). Let  us denote
$$
\mathcal{E}_1(f,f)
=  \int_D \int_D
\frac{(f(x) - f(y))^2}{|x - y|^{d + \alpha}} \varphi_1(x) \varphi_1(y) \, dx \, dy.
$$
When $\mathcal{E}_1(f,f) = \infty$ then (\ref{Dirichlet0}) follows from (\ref{Dirichlet3}) by the Fatou lemma. Now let us consider the case $\mathcal{E}_1(f,f) < \infty$. By (\ref{transition1}) for any $t > 0$ we have
\begin{equation}
\label{Dirichlet4}
\frac{p_D(t,x,y)}{t} (f(x) - f(y))^2 \varphi_1(x) \varphi_1(y) 
\le \frac{c (f(x) - f(y))^2}{|x - y|^{d + \alpha}} \varphi_1(x) \varphi_1(y).
\end{equation}
The integral over $D \times D$ of the right-hand side of
 (\ref{Dirichlet4}) is
 equal to $c \mathcal{E}_1(f,f) < \infty$. Now (\ref{Dirichlet0})
 follows from (\ref{Dirichlet3}) by the bounded convergence theorem.

\end{proof}

\begin{proof}[Proof of Theorem \ref{variational}]
Let $f \in \mathcal{F}$. We have $f \varphi_1 \in L^2(D)$, 
$||f \varphi_1||_{L^2(D)} = 1$ and $f \varphi_1 \perp \varphi_1$ in
$L^2(D)$. Since $\{\varphi_n\}_{n = 1}^{\infty}$ is an orthonormal
basis in $L^2(D)$ we have
$$
f \varphi_1 = \sum_{n = 2}^{\infty} c_n \varphi_n,
$$
where $c_n = \int_D f(x) \varphi_1(x) \varphi_n(x) \, dx$ and the
equality holds in $L^2(D)$ sense. Hence 
$$
f = \sum_{n = 2}^{\infty} c_n \frac{\varphi_n}{\varphi_1}
$$
in $L^2(D,\varphi_1^2)$ sense. The condition $||f \varphi_1||_{L^2(D)}
= 1$ gives $\sum_{n = 1}^{\infty} c_n^2 = 1$.

We will show that
\begin{equation}
\label{variational2}
\mathcal{E}(f,f) = \sum_{n = 2}^{\infty} (\lambda_n - \lambda_1) c_n^2.
\end{equation}

 We have
$$%\begin{eqnarray*}
T_t^D \frac{\varphi_n}{\varphi_1}(x) =
\int_D \frac{e^{\lambda_1 t} p_D(t,x,y)}{\varphi_1(x) \varphi_1(y)}
 \frac{\varphi_n(y)}{\varphi_1(y)}
 \varphi_1^2(y)
\, dy 
=  e^{-(\lambda_n - \lambda_1) t} 
 \frac{\varphi_n(x)}{\varphi_1(x)}.
$$ %\end{eqnarray*}
Hence by Parseval formula
$$
(T_t^D f,f)_{L^2(D,\varphi_1^2)} =
  \sum_{n = 2}^\infty c_n^2  e^{-(\lambda_n - \lambda_1) t},
$$
%It follows that
%\begin{eqnarray*}
%(f - T_t^D f,f)_{L^2(D,\varphi_1^2)} &=&
%\sum_{n = 2}^{\infty} c_n^2  (1 - e^{-(\lambda_n - \lambda_1) t}).
%\end{eqnarray*}
so 
\begin{equation}
\label{variational3}
\mathcal{E}(f,f) = \lim_{t\to 0^+} (f - T_t^D f,f)_{L^2(D,\varphi_1^2)} =
\lim_{t \to 0^+} \sum_{n = 2}^{\infty} c_n^2 
\frac{1 - e^{-(\lambda_n - \lambda_1) t}}{t}.
\end{equation}

%Let $f_k = \sum_{n = 2}^{k} c_n \varphi_n/\varphi_1$. We have
%\begin{eqnarray*}
%T_t^D f_k(x) &=&
%\int_D \frac{e^{\lambda_1 t} p_D(t,x,y)}{\varphi_1(x) \varphi_1(y)}
%\sum_{n = 2}^{k} c_n \frac{\varphi_n(y)}{\varphi_1(y)} \varphi_1^2(y)
%\, dy \\
%&=& \frac{e^{\lambda_1 t}}{\varphi_1(x)} 
%\sum_{n = 2}^{k} c_n \int_D p_D(t,x,y) \varphi_n(y) \, dy \\
%&=& \sum_{n = 2}^{k} c_n e^{-(\lambda_n - \lambda_1) t} 
% \frac{\varphi_n(x)}{\varphi_1(x)}.
%\end{eqnarray*}
%Hence
%\begin{eqnarray*}
%&& (T_t^D f_k,f_k)_{L^2(D,\varphi_1^2)} =
%\int_D T_t^D f_k(x) f_k(x) \varphi_1^2(x) \, dx \\
%&=& \sum_{n = 2}^{k} \sum_{m = 2}^{k} c_n c_m e^{-(\lambda_n -
%  \lambda_1) t} \int_D \varphi_n(x) \varphi_m(x) \, dx \\
%&=&  \sum_{n = 2}^{k} c_n^2  e^{-(\lambda_n - \lambda_1) t}.
%\end{eqnarray*}
%So we obtain
%$$
%(f_k - T_t^D f_k,f_k)_{L^2(D,\varphi_1^2)}
%= \sum_{n = 2}^{k} c_n^2  (1 - e^{-(\lambda_n - \lambda_1) t}).
%$$
%It follows that
%\begin{eqnarray*}
%(f - T_t^D f,f)_{L^2(D,\varphi_1^2)} &=&
%\lim_{k \to \infty} (f_k - T_t^D f_k,f_k)_{L^2(D,\varphi_1^2)} \\
%&=& \sum_{n = 2}^{\infty} c_n^2  (1 - e^{-(\lambda_n - \lambda_1) t}).
%\end{eqnarray*}
%So 
%\begin{equation}
%\label{variational3}
%\mathcal{E}(f,f) = 
%\lim_{t \to 0^+} \sum_{n = 2}^{\infty} c_n^2 
%\frac{1 - e^{-(\lambda_n - \lambda_1) t}}{t}.
%\end{equation}

To show (\ref{variational2}) we have to justify the change of the
limit and the sum in (\ref{variational3}).
%bd zmienilem na moim zdaniem prostszy argument (ponizej w komentarzach oryginalny)
%tk nie rozumiem dlaczego 
%tk e^{-(\lambda_n - \lambda_1) t})/t \uparrow \lambda_n - \lambda_1
%tk tzn rozumiem ze dla ustalonego n rosnie 
%tk ale chyba trzeba zeby istnialo takie t_0 zeby dla wszystkich n i dla wszytkich t< t_0 roslo 
%tk przywrocilem wiec poprzednia wersje ale nie jestem pewien 
%tk za duzo juz o tym nie myslalem
 Note that $(1 -
e^{-(\lambda_n - \lambda_1) t})/t \uparrow \lambda_n - \lambda_1$ when
$t \downarrow 0$ by convexity of the exponential function.
 Hence (\ref{variational2}) follows from
(\ref{variational3}) by the monotone convergence theorem.
% Note that $(1 -
%e^{-(\lambda_n - \lambda_1) t})/t \to \lambda_n - \lambda_1$ when $t
%\to 0^+$. When $\sum_{n = 2}^{\infty} (\lambda_n - \lambda_1) c_n^2 =
%\infty$ then (\ref{variational2}) follows from (\ref{variational3}) by
%the Fatou lemma. Now assume that
%$\sum_{n = 2}^{\infty} (\lambda_n - \lambda_1) c_n^2 < \infty$. Note
%that $0 < (1 - e^{-(\lambda_n - \lambda_1) t})/t \le \lambda_n -
%\lambda_1$ for $t > 0$. Hence (\ref{variational2}) follows from
%(\ref{variational3}) by the bounded convergence theorem.

By (\ref{variational2}) we get 
$$
\mathcal{E}(f,f) = \sum_{n = 2}^{\infty} (\lambda_n - \lambda_1) c_n^2
\ge (\lambda_2 - \lambda_1)  \sum_{n = 2}^{\infty} c_n^2
= \lambda_2 - \lambda_1.
$$
Now Lemma~\ref{Dirichlet} shows that the infimum in
(\ref{variational1}) is bigger or equal to $\lambda_2 -
\lambda_1$. When we put $f = \varphi_2/\varphi_1$ ($c_2 = 1$, $c_n =
0$ for $n \ge 3$) we obtain
$\mathcal{E}(\varphi_2/\varphi_1,\varphi_2/\varphi_1) = \lambda_2 -
\lambda_1$. This shows that the infimum in (\ref{variational1}) is
equal to $\lambda_2- \lambda_1$ and is
achieved for $f = \varphi_2/\varphi_1$.
\end{proof}

\section{Geometric and Analytic Properties of $\varphi_1$}

At first we recall the result which is already proven in \cite{BK3}, Theorem~2.1. (Theorem~2.1 in \cite{BK3} was formulated for $\alpha = 1$ (the Cauchy process) but the proof works for all $\alpha \in (0,2]$.)
 
\begin{thm}
\label{unimodality}
Let $D \subset \Rt$ be a bounded convex  domain which is
symmetric relative to both coordinate axes.
Then we have 
\begin{itemize}
\item[(i)] $\varphi_1$ is continuous and strictly positive in $D$.
\item[(ii)] $\varphi_1$ is symmetric in $D$ with respect to both coordinate axes. That is, 
$\varphi_1(x_1 ,-x_2) = \varphi_1(x_1, x_2)$ and  $ \varphi_1(-x_1, x_2) = 
\varphi_1(x_1, x_2)$.
\item[(iii)] $\varphi_1$ is unimodal in $D$ with respect to both 
coordinate axes. That is,  if we take any
$a_2 \in (-1,1)$ and $p(a_2)> 0$ such that $(p(a_2),a_2) \in \partial D$,
then the function $v(x_1)= \varphi_1(x_1, a_2)$ defined on $(-p(a_2), p(a_2))$ 
is non--decreasing on $(-p(a_2),0)$
and  non--increasing on $(0,p(a_2))$. 
Similarly, if we take any $a_1 \in (-L,L)$ and $r(a_1)>0$ such that 
$(a_1,r(a_1)) \in \partial D$,
then the function  $u(x_2)=  \varphi_1(a_1,x_2)$  defined on $(-r(a_1),r(a_1))$ is
non--decreasing on $(-r(a_1),0)$ and non--increasing on $(0,r(a_1))$.
\end{itemize}
\end{thm}

Next, we prove the Harnack inequality for $\varphi_1$. Such inequality is well known (see e.g. Theorem~6.1 in \cite{BB}). Our purpose here is to give a proof which will give an explicit constant. We adopt the method from \cite{BK3}. 

At first we need to recall some standard facts concerning stable processes.

By $P_{r,x}(z,y)$ we denote the
Poisson kernel for the ball $B(x,r) \subset \Rd$, $r > 0$ for the stable
process.  That is,   
$$P^z(X(\tau_{B(x,r)}) \in A)=\int_A
P_{r,x}(z,y) \, dy,$$ where $z \in B(x,r)$, $A \subset B^c(x,r)$. We
have \cite{BGR}
\begin{equation}
\label{Poissonf}
P_{r,x}(z,y) = C_{\alpha}^{d} 
\frac{(r^2 - |z - x|^2)^{\alpha/2}}{(|y - x|^2 - r^2)^{\alpha/2} |y - z|^d},
\end{equation}
where $C_{\alpha}^{d} = \Gamma(d/2) \pi^{-d/2 - 1} \sin(\pi \alpha /2)$, $ z \in B(x,r)$ and $y \in \interior(B^c(x,r))$. 

%tk zmiana
It is well known (\cite{G} cf. also \cite{BK} formula (2.10)) that 
\begin{equation}
\label{EytauB}
E^y(\tau_{B(0,r)}) = C_{\alpha}^{d} (\A_{d,-\alpha})^{-1} (r^2 - |y|^2)^{\alpha/2},
\end{equation}
where $r > 0$ and $\A_{d,-\alpha}$ is given by (\ref{Adalpha}).

When $d > \alpha$ by $G_D(x,y) = \int_0^{\infty} p_D(t,x,y) \, dt$ we denote the Green function for the domain $D \subset \Rd$, $x,y \in D$. We have $G_{D}(x,y) < \infty$ for $x \ne y$. (For $d \le \alpha$ the Green function may be defined by a different formula but we will not use it in this paper).
 
It is well known (see \cite{BGR}) that
\begin{equation}
\label{GreenB(0,1)}
G_{B(0,1)}(z,y)= \frac{R_{d,\alpha}}{|z - y|^{d - \alpha}}
\int_{0}^{w(z,y)} \frac{r^{\alpha/2 -1} dr}{(r + 1)^{d/2}},
\quad z,y \in B(0,1),
\end{equation}
where 
$$
w(z,y) = (1 - |z|^2)(1 - |y|^2)/|z - y|^2
$$
and $R_{d,\alpha} =
\Gamma(d/2)/(2^\alpha \pi^{d/2} (\Gamma(\alpha/2))^2)$.

By $\lambda_1(B_1)$ we denote the first eigenvalue for the unit ball $B(0,1)$.
Theorem~4 in \cite{BLM} (cf. also \cite{CS3}) gives the following estimate of $\lambda_1(B_1)$
\begin{equation}
\label{eigenvalue}
\lambda_1(B_1) \le (\mu_1(B_1))^{\alpha/2},
\end{equation}
where $\mu_1(B_1) \simeq 5.784$ is the first eigenvalue of the Dirichlet Laplacian for the unit ball. 

We will also need the following easy scaling property of $\varphi_1$.
\begin{lem}
\label{scalingphi1}
Let $D \subset \Rd$ be a bounded domain, $s > 0$ and $\varphi_{1,s}$ the first eigenfunction on the set $sD$ for the stable semigroup $\{P_{t}^{sD}\}_{t \ge 0}$. Then for any $x \in D$ we have
$\varphi_{1,s}(sx) = s^{-d/2} \varphi_{1,1}(x)$.
\end{lem}

Now we can formulate the Harnack inequality for $\varphi_1$.
\begin{thm}
\label{Harnack}
Let $\alpha \in (0,2)$, $d > \alpha$ and $D \subset \Rd$ be a bounded domain with inradius $R > 0$ and $0 < a < b < 1$. If $B(x,bR) \subset D$ then on $B(x,aR)$ $\varphi_1$ satisfies the Harnack inequality with constant $C_1 = C_1(d,\alpha,a,b)$. 
That is, for any $z_1, z_2 \in
B(x,aR)$ we have $\varphi_1(z_1) \le C_1 \varphi_1(z_2)$ where 
$$
C_1 = \frac{(b + a)^{d - \alpha/2} b^{\alpha}}{(b - a)^{d + \alpha/2}}
\left(1 + e + \frac{b^{d + \alpha/2} C_2}{(b - a)^{\alpha/2} (1 - b^{\alpha})^{d/\alpha}}\right)
$$ 
and $C_2 = C_2(d,\alpha) = \alpha^2 2^{3d/2 - \alpha/2 -1} C_{\alpha}^{d} M_{d,\alpha} (\lambda_1(B_1))^{d/\alpha}/((d - \alpha) R_{d,\alpha} \A_{d,-\alpha})$.
\end{thm}

\begin{proof}[Proof of Theorem \ref{Harnack}]

In view of Lemma~\ref{scalingphi1} we may and do assume that $R = 1$.

Let $B \subset D$ be any ball ($B \ne D$). For any $x,y \in B$, $t > 0$ we have
\begin{equation}
\label{expansion}
p_B(t,x,y) = \sum_{n = 1}^{\infty} e^{- \lambda_n(B) t} 
\varphi_{n, B}(x) \varphi_{n, B}(y), 
\end{equation}
where $\lambda_n(B)$ and $\varphi_{n, B}$ are the eigenvalues and
eigenfunctions for the semigroup $\{P_{t}^{B}\}_{t \ge 0}$. 

We will use the fact that the first eigenfunction is $q$-harmonic
in $B$ according to the $\alpha$-stable Schr\"odinger operator.

Let $\varphi_1$, $\lambda_1 = \lambda_1(D)$ be the first eigenfunction 
and eigenvalue
for the semigroup $\{P_{t}^{D}\}_{t \ge 0}$. Let $A$ be the
infinitesimal generator of this semigroup. For $x \in D$ we have 
$$
A \varphi_1(x) = 
\lim_{t \to 0^+} \frac{P_t^D \varphi_1(x) - \varphi_1(x)}{t} = 
\frac{e^{- \lambda_1(D) t} \varphi_1(x) - \varphi_1(x)}{t} =
- \lambda_1(D) \varphi_1(x). 
$$
This gives that $(A + \lambda_1(D))\varphi_1 = 0$ on $D$. It follows that
$\varphi_1$ is $q$-harmonic on $B$ according to the $\alpha$-stable Schr\"odinger
operator $A + q$  with $q \equiv \lambda_1(D)$. Formally this follows
from Proposition~3.17, Theorem~5.5,
 Definition~5.1 from \cite{BB} and the fact that $(B,\lambda_1(D))$ is
 gaugeable because $B$ is a proper open subset of
 $D$ and $\lambda_1(B) > \lambda_1(D)$.

Let $V_B(x,y) = \int_0^{\infty} e^{\lambda_1(D) t} p_B(t,x,y) \, dt$.
Here, $V_B$ is the  $q$-Green function, for $q \equiv \lambda_1(D)$, see page 58 in \cite{BB}. 
The $q$-harmonicity of $\varphi_1$ (Definition 5.1 in \cite{BB}), Theorem~4.10
 in \cite{BB} (formula (4.15)) and
formula (2.17) in \cite{BB} (page 61) give that for $z \in B$, 
\begin{eqnarray}
\nonumber
\varphi_1(z) 
&=&
E^z[e_{\lambda_1(D)}(\tau_{B}) \varphi_1(X(\tau_{B}))] \\
\label{qharmonic}
&=&
\A_{d,-\alpha} \int_B V_B(z,y) \int_{D \setminus B} |y - w|^{-d - \alpha} 
\varphi_1(w) \, dw \, dy, 
\end{eqnarray}   
where $e_{\lambda_1(D)}(\tau_{B}) = \exp(\lambda_1(D) \tau_B)$. 
Of course (\ref{qharmonic}) is a standard fact in the theory of
$q$-harmonic functions  for the $\alpha$-stable Schr\"odinger
operators. 
For us this will be a  key formula for proving the  
Harnack inequality for $\varphi_1$. 

By the well  known formula for the distribution of the harmonic
measure \cite{IW} we have 
\begin{equation}
\label{Ikeda}
E^{z}\varphi_1(X(\tau_B)) = 
\A_{d,-\alpha} \int_B G_B(z,y) \int_{D \setminus B} |y - w|^{-d - \alpha} 
\varphi_1(w) \, dw \, dy.
\end{equation} 

To obtain our Harnack inequality for $\varphi_1$ we will  first compare 
(\ref{qharmonic}) and (\ref{Ikeda}) and then we will use the formula for 
$E^{z}\varphi_1(X(\tau_B))$. In order to compare 
(\ref{qharmonic}) and (\ref{Ikeda}) we need to compare 
$V_B(z,y)$ and $G_B(z,y)$. This will be
done in a sequence of lemmas.

\begin{lem}
\label{VBzy1}
Let $D \subset \Rd$, $d > \alpha$ be a bounded domain with inradius $1$ and $B
\subsetneq D$ be a ball with radius $b < 1$. 
Then for any $z,y \in B$ and $t_0 > 0$ we have
$$
V_B(z,y) \le e^{\lambda_1(B_1) t_0} \int_{0}^{t_0} p_B(t,z,y) \, dt +
\frac{C_3}{ t_0^{(d - \alpha)/\alpha}},
$$
where $B_1 = B(0,1)$ and $C_3 =
\alpha (d - \alpha)^{-1} (1 - b^{\alpha})^{-d/\alpha} M_{d,\alpha}$.
\end{lem}

\begin{proof}
The inradius of $D$ is $1$ so $\lambda_1(D) \le \lambda_1(B_1)$. It
follows that
\begin{equation}
\label{twointegrals}
V_B(z,y) \le e^{\lambda_1(B_1) t_0} \int_{0}^{t_0} p_B(t,z,y) \, dt +
\int_{t_0}^{\infty} e^{\lambda_1(B_1)t} p_B(t,z,y) \, dt.
\end{equation}
By (\ref{expansion}) we obtain
$$
p_B(t,z,y) = \sum_{n =1}^{\infty} e^{-\lambda_n(B)t} \varphi_{n,B}(z) \varphi_{n,B}(y)
\le \frac{1}{2} \sum_{n = 1}^{\infty} e^{-\lambda_n(B)t} (\varphi_{n,B}^2(z) + \varphi_{n,B}^2(y)).
$$
It follows that the second integral in (\ref{twointegrals}) is bounded  above by
\begin{equation}
\label{expansionBD}
\frac{1}{2} \int_{t_0}^{\infty} \sum_{n=1}^{\infty}
e^{(\lambda_1(B_1) - \beta \lambda_n(B))t} e^{- \lambda_n(B)(1-\beta) t}
(\varphi_{n,B}^2(z) + \varphi_{n,B}^2(y)) \, dt,
\end{equation}
where $\beta = \lambda_1(B_1)/\lambda_1(B) = b^{\alpha}$ (see \ref{gapscaling}). 

Note also that
$e^{\lambda_1(B_1) - \beta \lambda_n(B)} \le e^{\lambda_1(B_1) - \beta
  \lambda_1(B)} = e^{0} = 1$.

For any $w \in B$ ($w = z$ or $w = y$) we  have
\begin{eqnarray*}
&&
\int_{t_0}^{\infty} \sum_{n =0}^{\infty} e^{-\lambda_n(B) (1-\beta) t} \varphi_{n,B}^2(w) \, dt =
\int_{t_0}^{\infty} p_B((1 - \beta)t,w,w) \, dt \\
&&
\le \int_{t_0}^{\infty} p((1 - \beta)t,0,0) \, dt  
\le \int_{t_0}^{\infty} \frac{M_{d,\alpha}}{(1 -\beta)^{d/\alpha}
  t^{d/\alpha}} \, dt
= \frac{C_3}{t_0^{(d - \alpha)/\alpha}}.
\end{eqnarray*}
\end{proof}

\begin{lem}
\label{EG}
Let $0<a<b<1$, $B =  B(w,b)$, $w \in \Rd$. For any $y \in B$
and $z \in B(w,a)$ we have
$$
C_4 G_B(z,y) \ge  E^{y}(\tau_B),
$$
where 
$C_4 = b^{d + \alpha/2} \alpha 2^{3d/2 - \alpha/2 - 1} C_{\alpha}^{d}/((b - a)^{\alpha/2} R_{d,\alpha} \A_{d,-\alpha})$.
\end{lem}
\begin{proof}
We may and do assume that $w = 0$. Let us consider the formula for the Green function for a unit ball $G_{B(0,1)}(z,y)$ (\ref{GreenB(0,1)}).
Note that for any $t >  0$ 
$$
\int_{0}^{t} \frac{r^{\alpha/2 -1} dr}{(r + 1)^{d/2}}
\ge \frac{1}{2^{d/2}} \int_{0}^{t \wedge 1} r^{\alpha/2 - 1}
= \frac{ (t^{\alpha/2} \wedge 1)}{\alpha 2^{d/2 - 1}}.
$$
Hence for any $z,y \in B(0,1)$
$$
G_{B(0,1)}(z,y) \ge R_{d,\alpha} \alpha^{-1} 2^{-d/2 + 1} |z -
y|^{\alpha - d} (1 \wedge (w(z,y))^{\alpha/2}). 
$$

By scaling it follows that for any $z,y \in B$,
\begin{eqnarray}
\nonumber
G_B(z,y) &=&
 b^{\alpha - d} G_{B(0,1)} \left(\frac{z}{b},\frac{y}{b}\right) \\
\nonumber
&\ge& 
\frac{R_{d,\alpha} \alpha^{-1} 2^{-d/2 + 1}}{b^{d - \alpha}
  \left|\frac{z}{b} - \frac{y}{b}\right|^{d - \alpha}}
\left(1 \wedge \frac{\left(1 - \left|\frac{z}{b}\right|^2 \right)^{\alpha/2} 
\left(1 - \left|\frac{y}{b}\right|^2 \right)^{\alpha/2}}{\left|\frac{z}{b}
  - \frac{y}{b}\right|^{\alpha}}\right)
\\
\label{Greenest}
&=& \frac{R_{d,\alpha} \alpha^{-1} 2^{-d/2 + 1}}{b^{\alpha} |z -
  y|^{d - \alpha}}
\left(b^{\alpha} \wedge \frac{(b^2 - |z|^2)^{\alpha/2}(b^2 -
  |y|^2)^{\alpha/2}}{|z - y|^{\alpha}}\right).
\end{eqnarray}
For $z \in B(0,a)$ and $y \in B = B(0,b)$ we have $|z - y| \le a + b
\le 2b$ and $(b^2 - |z|^2)^{\alpha/2} \ge (b^2 - a^2)^{\alpha/2}$.
Hence
$$
\frac{(b^2 - |z|^2)^{\alpha/2}}{|z - y|^{\alpha}}
\ge \frac{((b - a)(b + a))^{\alpha/2}}{((a + b)^2)^{\alpha/2}}
\ge \frac{1}{2^{\alpha/2}} \left(1 - \frac{a}{b}\right)^{\alpha/2}.
$$
It follows that for $z \in B(0,a)$ and $y \in B(0,b)$,
(\ref{Greenest}) is bounded below by
$$
\frac{R_{d,\alpha} \alpha^{-1} 2^{-d/2 + 1}}{b^{d}2^{d - \alpha}
  2^{\alpha/2}}   \left(1 - \frac{a}{b}\right)^{\alpha/2} (b^2 -
  |y|^2)^{\alpha/2}.
$$
By the formula for $E^{y}(\tau_B)$ (\ref{EytauB}) this is equal to $C_4^{-1} E^{y}(\tau_B)$. 
\end{proof}

\begin{lem}
\label{VBzy2}
Let $D \subset \Rd$, $d > \alpha$ be a bounded domain with inradius $1$, $0< a < b < 1$ and $B = B(x,b) \subset D$. Then for any $z \in
B(x,a)$ and $y \in B$ we have $G_B(z,y) \le V_{B}(z,y) \le C_5 \,
G_B(z,y)$,where $C_5 = 1 + e + C_3 C_4 (\lambda_1(B_1))^{d/\alpha}$.
\end{lem}
\begin{proof}
The inequality $G_B(z,y) \le V_B(z,y)$ is trivial, it follows from the
definition of $G_B(z,y)$ and $V_B(z,y)$. 

We will prove the inequality $ V_{B}(z,y) \le C_5 \, G_B(z,y)$.
By Lemma~4.8 in \cite{BB} we have
\begin{equation}
\label{recurrence}
V_{B}(z,y) = G_{B}(z,y) + \lambda_1(D) \int_{B} V_{B}(z,u) G_{B}(u,y) \, du.
\end{equation}
By Lemma~\ref{VBzy1},  $\int_B V_{B}(z,u) G_B(u,y) \, du$ is bounded  above by
\begin{equation}
\label{VGest}
e^{\lambda_1(B_1) t_0} \int_{B} \int_{0}^{t_0} p_{B}(t,z,u) \, dt
G_{B}(u,y) \, du +
\frac{C_3}{t_0^{(d - \alpha)/\alpha}} \int_B
G_B(u,y) \, du.
\end{equation}
Let us denote the above sum
by $\text{I} + \text{II}$. We have
\begin{eqnarray*}
&&
\int_B \int_{0}^{t_0} p_B(t,z,u) \, dt \, G_B(u,y) \, du 
= \int_{0}^{t_0} \int_{0}^{\infty} \int_{B} p_{B}(t,z,u) p_{B}(s,u,y)
\, du \, ds \, dt
\\
&&
= \int_{0}^{t_0} \int_{0}^{\infty} p_B(t + s,z,y) \, ds \, dt
\le t_0 G_B(z,y).
\end{eqnarray*}
It follows that
$\text{I} \le t_0 e^{\lambda_1(B_1) t_0} G_B(z,y)$. 

By applying Lemma~\ref{EG} for $z \in B(x,a)$ we get
$$
\text{II} = 
\frac{C_3 E^{y}(\tau_B)}{t_0^{(d - \alpha)/\alpha}}
\le
\frac{C_3 C_4 G_B(z,y)}{t_0^{(d - \alpha)/\alpha}}
$$
Putting the estimates (\ref{recurrence}), (\ref{VGest}) together with those
for   $\text{I}$ and $\text{II}$  gives 
\begin{equation}
\label{finalV}
V_B(z,y) \le G_B(z,y) \left(1 + \lambda_1(B_1) t_0 e^{\lambda_1(B_1)
    t_0} + 
\frac{C_3 C_4 \lambda_1(B_1)}{ t_0^{(d - \alpha)/\alpha}}\right).
\end{equation}

Putting $t_0 = 1/\lambda_1(B_1)$ we obtain
$$
V_B(z,y) \le G_B(z,y) (1 + e + C_3 C_4 (\lambda_1(B_1))^{d/\alpha}).
$$

\end{proof}

We now return  to the proof of Theorem~\ref{Harnack}.
Let $z_1, z_2 \in B(x,a) \subset B(x,b) \subset D$. By
(\ref{qharmonic}), (\ref{Ikeda}) and Lemma~\ref{VBzy2} we obtain
\begin{equation}
\label{phiz2}
\varphi_1(z_2) \ge E^{z_2}[\varphi_1(X(\tau_{B(x,b)}))]
\end{equation}
and
\begin{equation}
\label{phiz1}
\varphi_1(z_1) \le C_5 \, E^{z_1}[\varphi_1(X(\tau_{B(x,b)}))].
\end{equation}
So to compare $\varphi_1(z_2)$ and $\varphi_1(z_1)$ we have to compare 
$E^{z_1}[\varphi_1(X(\tau_{B(x,b)}))]$ and 
$E^{z_2}[\varphi_1(X(\tau_{B(x,b)}))]$. 

We have 
\begin{equation}
\label{Poisson}
E^{z_i}[\varphi_1(X(\tau_{B(x,b)}))] = 
\int_{D \setminus \overline{B(x,b)}} \varphi_1(y) P_{b,x}(z_i,y) \, dy,
\end{equation}
for $i = 1,2$, where $P_{b,x}(z_i,y)$ is the Poisson kernel for the ball $B(x,b)$ which is given by an explicit formula (\ref{Poissonf}). 
We have reduce  to comparing  $P_{b,x}(z_1,y)$ and
$P_{b,x}(z_2,y)$. Recall that $z_1, z_2 \in B(x,a)$. For $y \in
B^{c}(x,b)$ we have
$$
\frac{|y-z_2|}{|y-z_1|} \le \frac{b + a}{b - a} 
$$
and
$$
\frac{(b^2 - |z_1 - x|^2)^{\alpha/2}}{(b^2 - |z_2 - x|^2)^{\alpha/2}} 
\le 
\frac{b^{\alpha}}{(b^2 - a^2)^{\alpha/2}}.
$$
It follows that 
$$
\frac{P_{b,x}(z_1,y)}{P_{b,x}(z_2,y)} \le
\frac{(b + a)^{d - \alpha/2} b^{\alpha}}{(b - a)^{d + \alpha/2}}.
$$
Using this, (\ref{Poisson}), (\ref{phiz1}) and (\ref{phiz2}) we obtain
for $z_1, z_2 \in B(x,a)$  
$$
\varphi_1(z_1) \le C_5 (b + a)^{d - \alpha/2} b^{\alpha} (b - a)^{-d - \alpha/2} \varphi_1(z_2).
$$
\end{proof}

In this paper we will need the Harnack inequality for $\varphi_1^2$ in
dimension $d = 2$. For this reason we will formulate the following
corollary of Theorem~\ref{Harnack}. In this corollary we choose $b \in
(0,1/2]$ and $a = b/2$.

\begin{cor}
\label{corHarnack}
Let $\alpha \in (0,2)$ and $D \subset \Rt$ be a bounded domain with
inradius $R > 0$ and $b \in (0,1/2]$. If $B(x,bR) \subset D$ then on
$B(x,bR/2)$ $\varphi_1^2$
satisfies the Harnack inequality with constant $c_H =
c_H(\alpha)$. 
That is, for any $z_1, z_2 \in
B(x,bR/2)$ we have $\varphi_1^2(z_1) \le c_H \varphi_1^2(z_2)$ where 
\begin{equation}
\label{corHarnack1}
c_H = 3^{4 - \alpha} 2^{2 \alpha} \left(4 + \frac{12 \Gamma(2/\alpha)}{\alpha (2 - \alpha) (1 - 2^{-\alpha})^{2/\alpha}}\right)^2. 
\end{equation}
We point out that $c_H$ does not depend on $b \in (0,1/2]$.
\end{cor}
\begin{proof}
We are going to obtain upper bound estimates for constants $C_1$, $C_2$ from Theorem~\ref{Harnack} for $d = 2$, $a = b/2$ and $b \in (0,1/2]$. 

Putting d = 2 we get $C_{\alpha}^2 = \pi^{-2} \sin(\pi \alpha/2)$, $R_{2,\alpha} = 2^{-\alpha} \pi^{-1} \Gamma^{-2}(\alpha/2)$, $M_{2,\alpha} = 2^{-1} \pi^{-1} \alpha^{-1} \Gamma(2/\alpha)$, $\mathcal{A}_{2,-\alpha} = \alpha^2 2^{\alpha - 2} \pi^{-1} \Gamma(\alpha/2) \Gamma^{-1}(1 - \alpha/2)$. 

Putting these constants to the formula for $C_2$ and using also the fact that $\Gamma(\alpha/2) \Gamma(1 - \alpha/2) = \pi \sin^{-1}(\pi \alpha/2)$ we obtain after easy calculations
$$
C_2 = \frac{2^{3 - \alpha/2} \Gamma(2/\alpha) (\lambda_1(B_1))^{2/\alpha}}{(2 - \alpha) \alpha} 
\le \frac{6 \cdot 2^{3 - \alpha/2} \Gamma(2/\alpha)}{(2 - \alpha) \alpha}.
$$
The last inequality follows from (\ref{eigenvalue}) and the fact that $\mu_1(B_1) < 6$. 

Putting $d = 2$ and $a = b/2$ we obtain
$$
C_1 = 3^{2 - \alpha/2} 2^{\alpha} \left(1 + e + \frac{2^{\alpha/2} b^2 C_2}{(1 - b^{\alpha})^{2/\alpha}}\right).
$$
Now using  the estimate for $C_2$ and the inequality $b \le 1/2$ we get
\begin{equation}
\label{corHarnack2}
C_1 \le 3^{2 - \alpha/2} 2^{\alpha} \left(4 + \frac{12 \Gamma(2/\alpha)}{\alpha (2 - \alpha) (1 - 2^{-\alpha})^{2/\alpha}}\right). 
\end{equation}
In the assertion of Corollary \ref{corHarnack} we have the Harnack inequality for $\varphi_1^2$ so $c_H$ is equal to the square of the right hand side of (\ref{corHarnack2}).
\end{proof}

\section{Spectral gap for rectangles} \label{sec:above}

We begin from several lemmas, which will lead us to the estimation
of the spectral gap for rectangles.
\begin{lem}\label{phi1est}
Let $D=(-L,L)\times(-1,1)$, where $L\geq 1$. Then
$$
  \varphi_1(x) \leq \frac{3}{\sqrt{L}} \quad\textrm{for all $x\in D$}
$$
and
$$
  \varphi_1(x_1,x_2) \geq \frac{1}{2 \sqrt{L}}
   (1-\frac{2}{L}|x_1|)(1-2|x_2|)
$$
 for all $(x_1,x_2)\in [-L/2,L/2]\times [-1/2,1/2]$.
\end{lem}
\begin{proof}
The lemma easily follows from unimodality and symmetry of $\varphi_1$ (see Theorem~\ref{unimodality}), midconcavity of $\varphi_1$ (see Theorem~1.1 in \cite{BKM}) and the equality $\int_D \varphi_1^2\, dx=1$.
\end{proof}

\begin{lem}\label{partition0}
Let $\mu_k>0$ ($k=1, \ldots, L$),
%tk dodalem L \ge 2
$L \ge 2$
 be unimodal, i.e.,
there exists $k_0$ such that $\mu_i \leq \mu_j$ for $i\leq j\leq k_0$
and $\mu_i \geq \mu_j$ for $k_0\leq i\leq j$.
Then for any $f_k\in\R$ such that $\sum_{k=1}^L f_k\mu_k=0$ we have
$$
 \sum_{k=1}^L \mu_k f_k^2 \leq L^2
   \sum_{k=1}^{L - 1} (\mu_k \wedge \mu_{k+1})(f_k-f_{k+1})^2.
$$
\end{lem}
\begin{proof}
Let $M=\sum_{k=1}^L \mu_k$.
By  $\sum_{k=1}^L \mu_k f_k=0$ and Schwarz inequality we obtain
\begin{eqnarray}
 M \sum_{k=1}^L \mu_k f_k^2
 &=& \sum_{j=1}^L \mu_j \sum_{k=1}^L \mu_k f_k^2 =
\frac{1}{2} \sum_{j,k=1}^L \mu_j \mu_k (f_j^2 + f_k^2) \nonumber\\
&=&
\frac{1}{2} \sum_{j,k=1}^L \mu_j \mu_k (f_j - f_k)^2 \nonumber\\
&=&
 \sum_{1 \leq j<k\leq L} \mu_j \mu_k \expr{\sum_{t=j}^{k-1} (f_t - f_{t+1}) }^2 \nonumber\\
&\leq&
L \sum_{1 \leq j<k\leq L} \mu_j \mu_k \sum_{t=j}^{k-1} (f_t - f_{t+1})^2 \nonumber\\
&=&
L \sum_{t=1}^{L-1}
\expr{ \sum_{j \leq t < k} \mu_j \mu_k } \cdot
 (f_t-f_{t+1})^2.  \label{juzprawie}
\end{eqnarray}

For $t < k_0$ (where $k_0$ is defined in the lemma) we have
$$
\sum_{j \leq t < k} \mu_j \mu_k \leq
L \mu_t \sum_{k=1}^L \mu_k =
L (\mu_t \wedge \mu_{t+1}) M.
$$
Similarly for $t \geq k_0$
$$
\sum_{j \leq t < k} \mu_j \mu_k \leq
L \mu_{t+1} \sum_{j=1}^L \mu_j =
L (\mu_t \wedge \mu_{t+1}) M.
$$
These two inequalities combined with (\ref{juzprawie}) finish the proof.
\end{proof}

\begin{lem}\label{partition}
Let $(D,\mu)$ be a finite measure space and
$D=\bigcup_{k=1}^L D_k$,
%tk dodalem L \ge 1
$L \ge 1$
 with pairwise disjoint $D_k$'s.
We assume that the sequence $\mu_k=\mu(D_k)>0$ is unimodal.
%, i.e.,
%there exists $k_0$ such that $\mu_i \leq \mu_j$ for $i\leq j\leq k_0$
%and $\mu_i \geq \mu_j$ for $k_0\leq i\leq j$.
Then
\begin{eqnarray}
%tk ponumerowalem wszystkie linie
&& \frac{1}{\mu(D)} \int_D\int_D (f(x)-f(y))^2 \mu(dx)\mu(dy) \label{partition1}\\
& \leq &
2 \sum_{k=1}^L
 \frac{1}{\mu_k} \int_{D_k}\int_{D_k} (f(x)-f(y))^2\mu(dx)\mu(dy) \label{partition2} 
\\
&&+
 4\, L^2 
\sum_{k=1}^{L-1}
 \frac{1}{\mu_k \vee \mu_{k+1}} \int_{D_k}\int_{D_{k+1}} (f(x)-f(y))^2\mu(dx)\mu(dy)
\qquad\label{partition3}
\end{eqnarray}
for all $f\in L^2(D,\mu)$.
\end{lem}

\begin{proof}
Let $f\in L^2(D,\mu)$.
Without loss of generality we may assume that
%tk L \ge 2
$L \ge 2$ and $\int_D f d\mu =0$.
Then 
%tkthe left hand side of (\ref{nierlematu}) 
(\ref{partition1}) is equal to 
$2  \int_D f^2 d\mu$.

Let $f_k=\frac{1}{\mu_k} \int_{D_k}f d\mu$.
We have
$$
\sum_{k=1}^L
 \frac{1}{\mu_k} \int_{D_k}\int_{D_k} (f(x)-f(y))^2\mu(dx)\mu(dy)
 = 2\int_D f^2 d\mu - 2 \sum_{k=1}^L \mu_k f_k^2.
$$
Thus if $\sum_{k=1}^L \mu_k f_k^2 \leq \frac{1}{2} \int_D f^2 d\mu$, then
%tk (\ref{nierlematu}) 
(\ref{partition1} - \ref{partition3}) holds.
Consequently, from now on we may assume that 
\begin{equation}\label{mikfk2}
\sum_{k=1}^L \mu_k f_k^2 > \frac{1}{2} \int_D f^2 d\mu.
\end{equation}
Thus, by Lemma~\ref{partition0} we have
\begin{equation}\label{zlematu}
2 \int_D f^2 d\mu <
4 L^2
\sum_{k=1}^{L-1} (\mu_k \wedge \mu_{k+1}) (f_k - f_{k+1})^2.
\end{equation}
On the other hand,
\begin{eqnarray*}
&& \int_{D_k}\int_{D_{k+1}} (f(x)-f(y))^2\mu(dx)\mu(dy) \\
 &=&
 \int_{D_k}\int_{D_{k+1}} ((f(x)-f_k)-(f(y)-f_{k+1}) + (f_k-f_{k+1}))^2\mu(dx)\mu(dy) \\
&\geq& \mu_k \mu_{k+1} (f_k - f_{k+1})^2.
\end{eqnarray*}
The lemma now follows from (\ref{zlematu}).
\end{proof}

\begin{lem}
\label{recHarnack}
Let $\alpha \in [1,2)$, $D = (-L,L) \times (-1,1)$, $L \ge 1$ and
$\varphi_1$ be the first eigenfunction for $\{P_t^D\}_{t \ge 0}$. Let $-L + 1/4 \le a \le b \le L - 1/4$, $b - a = 1/8$ and put $A = [a, b] \times
[-1/8,1/8]$. Then we  have
\begin{equation}
\label{recHarnackf1}
\left( \max_{x \in A} \varphi_1^2(x) \right)
\left(\int_A \varphi_1^2(x) \, dx \right)^{-1} \le C_R,
\end{equation}
where $C_R =  10^4$. 
\end{lem}
\begin{proof}
We will use the fact that $\varphi_1$ is symmetric and unimodal with respect to both coordinate axes (see Theorem~\ref{unimodality}). We will also use much stronger fact that $\varphi_1$ is "midconcave" (see Theorem~1.1 in \cite{BKM}). That is for any $x_2 \in (-1,1)$ $x_1 \to \varphi_1(x_1,x_2)$ is concave on $(-L/2,L/2)$ and for any $x_1 \in (-L,L)$ $x_2 \to \varphi_1(x_1,x_2)$ is concave on $(-1/2,1/2)$.

By symmetry of $\varphi_1$ we may and do assume that $b \ge 0$. We will consider two cases: Case 1, $b \in [0,3/8)$, Case 2, $b \in [3/8,L)$.

At first let us consider Case 1: $b \in [0,3/8)$.
Note that by the unimodality $\min_{x \in A} \varphi_1(x)$ is equal to $\varphi_1(b,1/8)$ or $\varphi_1(a,1/8)$. By concavity of $x_2 \to \varphi_1(b,x_2)$ on $(-1/2,1/2)$ we obtain
$$
\varphi_1(b,1/8) \ge (3/4) \varphi_1(b,0) + (1/4) \varphi_1(b,1/2) \ge (3/4) \varphi_1(b,0).
$$
On the other hand $x_1 \to \varphi_1(x_1,0)$ is concave on $(-L/2,L/2)$. We have $L \ge 1$ so $x_1 \to \varphi_1(x_1,0)$ is concave on $(-1/2,1/2)$. It follows that
$$
\varphi_1(b,0) \ge \varphi_1(3/8,0) \ge (1/4) \varphi_1(0,0) + (3/4) \varphi_1(1/2,0) \ge (1/4) \varphi_1(0,0).
$$
Hence $\varphi_1(b,1/8) \ge (3/16) \varphi_1(0,0)$. Similarly one can show that $\varphi_1(a,1/8) \ge (3/16) \varphi_1(0,0)$. 

We also have $\max_{x \in D} \varphi_1(x) = \varphi_1(0,0) \ge \max_{x \in A} \varphi_1(x)$. Finally 
$$
\int_{A} \varphi_1^2(x) \, dx \ge |A| \min_{x \in A} \varphi_1^2(x) \ge \frac{1}{32} \left(\frac{3}{16}\right)^2 \varphi_1^2(0,0) \ge \frac{9}{2^{13}} \max_{x \in A} \varphi_1^2(x).
$$
This gives (\ref{recHarnackf1}) and finishes Case 1.

Now let us consider Case 2: $b \in [3/8,L)$.
Note that $\max_{x \in A} \varphi_1(x) = \varphi_1(a,0)$ and $\min_{x \in A} \varphi_1(x) = \varphi_1(b,1/8)$. As before $\varphi_1(b,1/8) \ge (3/4) \varphi_1(b,0)$. Hence
\begin{equation}
\label{rHarnack1}
\int_{A} \varphi_1^2(x) \, dx \ge |A| \left( \frac{3}{4} \right)^2 \varphi_1^2(b,0) =
\frac{9}{2^9} \varphi_1^2(b,0).
\end{equation}
Now we have to estimate $\varphi_1(b,0)$.

Let $x_0 = (b,0)$, $r = \sqrt{2}/8$ and consider a ball $B = B(x_0,r)$. It is easy to note that $\overline{B} \subset D$. By formula (\ref{qharmonic}) and the fact that $e_{\lambda_1(D)}(\tau_{B}) \ge 1$ we have 
\begin{equation}
\label{rHarnack2}
\varphi_1(x_0) \ge 
E^{x_0}[ \varphi_1(X(\tau_{B}))].
\end{equation}

Now let us introduce polar coordinates $(\rho,\psi)$ with centre at $x_0 = (b,0)$. For any $z = (z_1,z_2) \in \Rt$ we have $z_1 - b = \rho \cos(\psi)$, $z_2 = \rho \sin(\psi)$. Let us consider the set $S_1 = \{(\rho,\psi): \rho \in (r,2 r), \, \psi \in (3 \pi/4, 5 \pi/4)\}$. Note that $S_1$ is chosen so that $S_1 \subset [b - 2 r,b - r/\sqrt{2}] \times [-\sqrt{2} r, \sqrt{2} r] \subset [b - 3/8,a] \times [-1/4,1/4]$. 

By unimodality and "midconcavity" for any $z \in [b - 3/8,a] \times [-1/4,1/4]$ we have $\varphi_1(z) \ge \varphi_1(a,0)/2$. This and (\ref{rHarnack2}) gives 
\begin{equation}
\label{rHarnack3}
\varphi_1(x_0) \ge 
E^{x_0}[ \varphi_1(X(\tau_{B})); \, X(\tau_B) \in S_1] \ge 
(\varphi_1(a,0)/2) P^{x_0}( X(\tau_B) \in S_1).
\end{equation}
We have
$$
P^{x_0}( X(\tau_B) \in S_1) = \int_{S_1} P_{r,x_0}(x_0,y) \, dy,
$$
where $P_{r,x_0}(x_0,y)$ is the Poisson kernel for $B$ given by (\ref{Poissonf}). 

Let $S_2 = \{(\rho,\psi): \rho \in (2 r, \infty), \, \psi \in (3 \pi/4, 5 \pi/4)\}$. Since $P^{x_0}(X(\tau_B) \in B^c) = 1$ and the distribution $P^{x_0}(X(\tau_B) \in \cdot)$ is invariant under rotation around $x_0$ it is easy to note that $P^{x_0}(X(\tau_B) \in S_1 \cup S_2) = 1/4$. Hence
\begin{equation}
\label{rHarnack4}
P^{x_0}( X(\tau_B) \in S_1) 
= \frac{1}{4} - P^{x_0}( X(\tau_B) \in S_2) 
= \frac{1}{4} - \int_{S_2} P_{r,x_0}(x_0,y) \, dy.
\end{equation}
We have
\begin{equation}
\label{rHarnack5}
\int_{S_2} P_{r,x_0}(x_0,y) \, dy = 
C_{\alpha}^2 \int_{3 \pi/4}^{5 \pi/4} \int_{2 r}^{\infty} \frac{r^{\alpha}}{(\rho^2 - r^2)^{\alpha/2} \rho^2} \rho \, d\rho \, d\psi.
\end{equation}
Note that $\rho^2 - r^2 \ge (3/4) \rho^2$ for $\rho \ge 2 r$ so (\ref{rHarnack5}) is smaller than
$$
C_{\alpha}^2 \frac{\pi}{2} r^{\alpha} \left(\frac{4}{3}\right)^{\alpha/2} \int_{2 r}^{\infty} \rho^{- \alpha - 1} \, d\rho 
= \frac{1}{2 \pi \alpha 3^{\alpha/2}} \sin\left(\frac{\pi \alpha}{2}\right) 
\le \frac{1}{2 \pi \sqrt{3}}.
$$
The last inequality follows from the fact that in this lemma we assume that $\alpha \in [1,2)$. 

Using this, (\ref{rHarnack3}) and (\ref{rHarnack4}) we obtain
$$
\varphi_1(x_0) \ge \frac{\varphi_1(a,0)}{2} \left(\frac{1}{4} - \frac{1}{2 \pi \sqrt{3}}\right).
$$
This and (\ref{rHarnack1}) gives
$$
\int_{A} \varphi_1^2(x) \, dx \ge 
\frac{9}{2^9} \left(\frac{1}{2} \left(\frac{1}{4} - \frac{1}{2 \pi \sqrt{3}}\right)\right)^2
\varphi_1^2(a,0) \ge \frac{\varphi_1^2(a,0)}{10^4} = C_R^{-1} \max_{x \in A} \varphi_1^2(x).
$$
\end{proof}

\begin{proof}[Proof of Theorem~\ref{estrectangle} -- part I]
By scaling of eigenvalues (see (\ref{gapscaling})) it is sufficient to show the following inequalities for rectangles $D = (-L,L) \times (-1,1)$, $L \ge 1$:
\begin{eqnarray}
\label{oszaczgory1}
  2\mathcal{A}_{2,-\alpha}^{-1} (
  \lambda_2-\lambda_1) \le 10^6 \cdot
  \left\{ \begin{array}{ll}
       \displaystyle
       \frac{2}{1-\alpha} \frac{1}{L^{1+\alpha}} & 
            \textrm{for $\alpha<1$,} \vspace{2mm}\\
    \displaystyle
      \frac{2 \log(L+1)}{L^2} & \textrm{for $\alpha=1$,}  \vspace{2mm}\\
\displaystyle
      (\frac{1}{2-\alpha}+\frac{1}{\alpha-1}) \frac{1}{L^2} & 
            \textrm{for $\alpha>1$.} 
\end{array} \right.
\end{eqnarray}
\begin{eqnarray}
\label{oszaczdolu1}
   2\mathcal{A}_{2,-\alpha}^{-1} (
  \lambda_2-\lambda_1) \geq 
\left\{ \begin{array}{ll}
  \displaystyle  
  \frac{1}{36\cdot 2^{1 + 2 \alpha}  (L)^{1+\alpha} }
  &\textrm{for $\alpha<1$,} \vspace{2mm}
  \\
  \displaystyle  
      10^{-9} \, \frac{\log(L+1)}{L^2} & \textrm{for $\alpha=1$,}  \vspace{2mm}\\
  \displaystyle
  \frac{1}{33 \cdot 13^{1 + \alpha/2} \cdot 10^4}\, \frac{1}{L^2}
  & \textrm{for $\alpha>1$.} 
    \end{array}
  \right.
\end{eqnarray}

Similarly, to prove Remark~\ref{remark1} it is sufficient to show
\begin{equation}
\label{remarkL}
2\mathcal{A}_{2,-\alpha}^{-1} (
  \lambda_2-\lambda_1) \geq 
  \frac{1}{36\cdot 2^{\alpha}  (L+1)^{1+\alpha} }.
\end{equation}

Let us take $f(x_1,x_2)=x_1$ for $x=(x_1,x_2)\in D$.
Then by Lemma~\ref{phi1est}
\begin{eqnarray*}
 \int_D f^2 \varphi_1^2 dx &\geq& \frac{1}{4L} \int_{-L/2}^{L/2} \int_{-1/2}^{1/2}
  (1-\frac{2}{L}|x_1|)^2 (1-2|x_2|)^2 x_1^2 \,dx_2 \,dx_1\\
  &=&
  \frac{L^2}{1440}.
\end{eqnarray*}
On the other hand, for $x\in D$ we have
\begin{eqnarray*}
&& \int_D |x-y|^{-\alpha} \,dy \leq
  \int_{B(0,\sqrt{5})} |y|^{-\alpha}\,dy\\
&& +
  \int_{-L}^L \int_{-1}^1 (|x_1-y_1|^2+|x_2-y_2|^2)^{-\alpha/2} 
    \ind_{\R\setminus[-1,1]}(y_1) \,dy_2 \, dy_1 \\
&& \leq
  2\pi\frac{5^{1-\alpha/2}}{2-\alpha} + 4\int_1^L y_1^{-\alpha}\,dy_1.
\end{eqnarray*}
Thus by Lemma~\ref{phi1est}
$$
 \mathcal{E}(f,f) <
  \frac{A_{2,-\alpha}}{2}
  \frac{9}{L} \cdot 4L \cdot (2\pi\frac{5^{1-\alpha/2}}{2-\alpha} +
   4\int_1^L y_1^{-\alpha}\,dy_1).
$$
Hence by Theorem~\ref{variational}
\begin{eqnarray*}
  \lambda_2-\lambda_1 &\leq& \frac{ \mathcal{E}(f,f)}{ \int_D f^2 \varphi_1^2 dx}\\
 &<& \frac{A_{2,-\alpha}}{2} \frac{72\cdot 1440}{L^2} 
     (\pi\frac{5^{1-\alpha/2}}{2-\alpha} + 2\int_1^L y_1^{-\alpha}\,dy_1),
\end{eqnarray*}
therefore (\ref{oszaczgory1}) is proven.

For $f\in L^2(D,\varphi_1^2)$ such that $\int_D f\varphi_1^2\, dx=0$ we have by
Lemma~\ref{phi1est}
\begin{eqnarray*}
   \int_D f^2 \varphi_1^2 \,dx &=&
\frac{1}{2} \int_D \int_D (f(x)-f(y))^2 \varphi_1^2(x)\varphi_1^2(y) \,dx\,dy \\
&\leq&
\frac{9 \diam(D)^{2+\alpha} }{2 L} 
   \int_D \int_D \frac{(f(x)-f(y))^2}{|x-y|^{2+\alpha}} \varphi_1(x)\varphi_1(y) \,dx\,dy \\
&\leq&
 36\cdot 2^{\alpha}  (L+1)^{1+\alpha} \frac{2}{A_{2,-\alpha}} \,
  \mathcal{E}(f,f),
\end{eqnarray*}
thus by Theorem~\ref{variational} we obtain 
(\ref{remarkL}) and also (\ref{oszaczdolu1})
in the case when $\alpha<1$.

Let $D=\bigcup_{k=1}^{[2L]} D_k$ be divided into $[2L]$ pairwise disjoint rectangles $D_k$
of size $\frac{2L}{[2L]} \times 2$,
denote $E_k=D_k\cup D_{k+1}$.
 Let $\mu = \varphi_1^2\,dx$, by Theorem~\ref{unimodality}
we see that $(D, \mu)$ satisfies the assumptions of Lemma~\ref{partition}.
Thus for $f\in L^2(D,\varphi_1^2)$ such that $\int_D f\varphi_1^2\, dx=0$ we have
\begin{eqnarray*}
%tk nie miescilo sie wiec zmienilem uklad
&&  2 \int_D f^2 \varphi_1^2 \,dx \leq
 2\sum_{k=1}^{[2L]}
 \frac{1}{\int_{D_k} \varphi_1^2\,dx} \int_{D_k}\int_{D_k} 
    (f(x)-f(y))^2\varphi_1^2(x) \varphi_1^2(y) \,dx\,dy \\
&&+ 4\, [2L]^2 
\sum_{k=1}^{[2L]-1}
%tk dodalem 2 w liczniku
%tk i pozniej dalej konsekwentnie zmienialem
 \frac{2}{\int_{E_k} \varphi_1^2\,dx} \int_{E_k}\int_{E_k} 
    (f(x)-f(y))^2\varphi_1^2(x) \varphi_1^2(y) \,dx\,dy \\
&& = I_1+I_2,
\end{eqnarray*}
\begin{eqnarray*}
I_2 &\leq&
%tk 32
 32 L^2
\sum_{k=1}^{[2L]-1}
 \frac{\sup_{E_k}\varphi_1^2}{\int_{E_k} \varphi_1^2\,dx} \\
&& \times
\int_{E_k}\int_{E_k} 
   (f(x)-f(y))^2
    \frac{\diam(E_k)^{2+\alpha}}{|x-y|^{2+\alpha}} \varphi_1(x) \varphi_1(y)
     \,dx\,dy. \label{quotient}
\end{eqnarray*}
%By unimodality and symmetry of $\varphi_1$, for each $k$ there exists
% $x_k\in E_k$ such that
%$B_k=B(x_k,1/4)\subset E_k$, $\sup_{B_k} \varphi_1^2 = \sup_{E_k} \varphi_1^2$,
%and $B(x_k,1)\subset D$. Thus by Corollary~\ref{Harnack1} we obtain
%$$
%\frac{\sup_{E_k}\varphi_1^2}{\int_{E_k} \varphi_1^2\,dx} \leq
% \frac{\sup_{B_k}\varphi_1^2}{\int_{B_k} \varphi_1^2\,dx} \leq
%   \frac{c_H}{|B_k|} = \frac{16\, c_H}{\pi}.
%$$
%Combining this and (\ref{quotient}) we get
Hence by Lemma~\ref{recHarnack} and $\diam(E_k)\leq \sqrt{13}$ we obtain
$$
   I_2 \leq
%tk zmiana 64 a nie 32 bo E_k razy E_k szacuja sie przez 2 razy E(f,f)
  64 \cdot 13^{1 + \alpha/2}C_R L^2 \,\frac{2}{A_{2,-\alpha}}\mathcal{E}(f,f).
$$
Similarly,
$$
 I_1 \leq
  2\cdot 13^{1+\alpha/2} C_R \,\frac{2}{A_{2,-\alpha}}\mathcal{E}(f,f)
$$
and (\ref{oszaczdolu1}) in the case when $\alpha>1$ follows.
\end{proof}

\begin{proof}[Proof of Theorem~\ref{estrectangle} -- part II, the case $\alpha=1$]
 Let $N = [L]$. We divide
$D$ into $2 N$  rectangles  of equal size $D_{-N+1},\ldots,D_{N}$, where 
$$
D_k = ((k - 1)L/N,kL/N) \times (-1,1), \quad k = -N + 1,\ldots, N.
$$
Let us note that Lemma~\ref{recHarnack} implies
$$
\left( \sup_{x \in D_k} \varphi_1^2 \right)
\left(\int_{D_k} \varphi_1^2 \right)^{-1} \le C_1^{-1},
$$
where $C_1 = C_R^{-1} = 10^{-4}$. 

We will also use the following easy inequality
$\inf_{x \in D_k, \, y \in D_{k + 1}} |x - y|^{-3} \ge C_2$, where $C_2 = (\sqrt{20})^{-3}$.

In the case $\alpha = 1$ we will show that
\begin{equation}
\label{reclnf0}
2 \A_{2,-1}^{-1} (\lambda_2 - \lambda_1) \ge \frac{C_1 C_2 \log(L + 1)}{360 \, L^2}.
\end{equation}
This implies (\ref{oszaczdolu1}) in the remaining case when $\alpha=1$. 

Fix $i \in \{1,\ldots,N\}$. For any $k = 1,\ldots, i$ let 
$$
A_k^i = \ldots \cup D_{k - 2i} \cup D_{k - i} \cup D_{k} \cup D_{k +
  i} \cup D_{k + 2 i} \cup \ldots
$$
Since $N$ is not necessarily divisible by $i$ the number of ``parts''
of $A_k^i$ may not be equal for different $k$. To make the definition
of $A_k^i$ more precise we introduce some more notation.

We have $N = i[N/i] + r(i)$ for some $r(i) \in \{0,\ldots,i-1\}$. 
Let $q(i,k) = [N/i]$ for $k = 1, \ldots, r(i)$, $q(i,k) = [N/i] -
1$ for $k = r(i) + 1,\ldots,i$ and $p(i,k) = -[N/i]$ for $k = 1,
\ldots,i - r(i)$, $p(i,k) = -[N/i] -
1$ for $k =i - r(i) + 1,\ldots,i$. 

For $m = p(i,k),\ldots,q(i,k)$ let $D_{k,m}^i = D_{k + m i}$. Then we
have
$$
A_k^i = \bigcup_{m = p(i,k)}^{q(i,k)} D_{k,m}^i. 
$$

Now we will apply Lemma~\ref{phi1est} to the set $A_k^i$ which is divided
as above. We take $\mu(dx) = \varphi_1^2(x) \, dx$ and $f =
\varphi_2/\varphi_1$. Let us denote $\mu_{k,m}^i = \int_{D_{k,m}^i}
\varphi_1^2$. 

Of course we have $\mu_{k,m}^i \vee \mu_{k,m+1}^i \ge 
(\mu_{k,m}^i  \mu_{k,m+1}^i)^{1/2}$ and
$$
 \int_{A_k^i}\int_{A_k^i} (f(x)-f(y))^2 \varphi_1^2(x) \varphi_1^2(y)
 \, dx \, dy
= 2 \int_{A_k^i} f^2 \varphi_1^2 \int_{A_k^i} \varphi_1^2 -
2 \left(\int_{A_k^i} f \varphi_1^2 \right)^2.
$$

So applying Lemma~\ref{phi1est} to $A_k^i$ and summing from $k = 1$ to $k = i$ we obtain
\begin{eqnarray}
&& 
\label{reclnf1}
2 \sum_{k = 1}^i   \int_{A_k^i} f^2 \varphi_1^2 
- 2 \sum_{k = 1}^i \left(\int_{A_k^i} \varphi_1^2 \right)^{-1}
  \left(\int_{A_k^i} f \varphi_1^2 \right)^2  \le \\
&& 
\label{reclnf2}
2  \sum_{k = 1}^i \sum_{m = p(i,k)}^{q(i,k)}
 \frac{1}{\mu_{k,m}^i} \int_{D_{k,m}^i}\int_{D_{k,m}^i} (f(x)-f(y))^2  
\varphi_1^2(x) \varphi_1^2(y) \, dx \, dy \\
&&
\label{reclnf3}
+
 4  \sum_{k = 1}^i \sum_{m = p(i,k)}^{q(i,k) - 1}
\frac{(q(i,k) - p(i,k) + 1)^2}{(\mu_{k,m}^i  \mu_{k,m+1}^i)^{1/2}}
 \int_{D_{k,m}^i} \int_{D_{k,m + 1}^i} \\
&&
\nonumber
\quad \quad \quad \quad \quad \quad \quad \quad\quad \quad \quad \quad
\quad \quad 
\times
(f(x)-f(y))^2 \varphi_1^2(x) \varphi_1^2(y) \, dx \, dy.
\end{eqnarray}

Now we will consider 2 cases:

Case 1. For any $i \in \{1,\ldots,[N^{1/4}]\}$ we have
\begin{equation}
\label{Case1}
\sum_{k = 1}^i \left(\int_{A_k^i} \varphi_1^2 \right)^{-1}
  \left(\int_{A_k^i} f \varphi_1^2 \right)^2 \le \frac{1}{2}.
\end{equation}

Case 2. There exists $i \in \{1,\ldots,[N^{1/4}]\}$ such that
\begin{equation}
\label{Case2}
\sum_{k = 1}^i \left(\int_{A_k^i} \varphi_1^2 \right)^{-1}
  \left(\int_{A_k^i} f \varphi_1^2 \right)^2 > \frac{1}{2}.
\end{equation}

At first we consider Case 1. Let us denote expressions in
(\ref{reclnf1}), (\ref{reclnf2}), (\ref{reclnf3}) by $L(i)$, $R(i)$,
$S(i)$ respectively. 

We have $\sum_{k = 1}^i \int_{A_k^i} f^2 \varphi_1^2 = \int_{D} f^2 \varphi_1^2 = 1$ so by the assumption (\ref{Case1}) we have $L(i) \ge 1$.

Now let us assume that for some $i \in \{1,\ldots,[N^{1/4}]\}$ we have $R(i) \ge S(i)$. This gives $R(i) \ge L(i)/2 \ge 1/2$. On the other hand we have
\begin{eqnarray}
\label{reclnf4}
&&
2 \A_{2,-1}^{-1} (\lambda_2 - \lambda_1) =
\int_{D} \int_{D}
\frac{(f(x) - f(y))^2}{|x - y|^3} \varphi_1(x) \varphi_1(y) \, dx \, dy
\\
\label{reclnf5}
&&
\ge \sum_{k = 1}^i \sum_{m = p(i,k)}^{q(i,k)} \int_{D_{k,m}^i} \int_{D_{k,m}^i}
\frac{(f(x) - f(y))^2}{|x - y|^3} \varphi_1(x) \varphi_1(y) \, dx \, dy. 
\end{eqnarray}
By our standard arguments (\ref{reclnf5}) is bounded below by
$$
C_1 C_2 \sum_{k = 1}^i \sum_{m = p(i,k)}^{q(i,k)} 
(\mu_{k,m}^i)^{-1} \int_{D_{k,m}^i} \int_{D_{k,m}^i}
(f(x) - f(y))^2  \varphi_1^2(x) \varphi_1^2(y) \, dx \, dy. 
$$
This is equal to $(C_1 C_2/2) R(i)$ where $R(i)$ is the expression in (\ref{reclnf2}). Since $R(i) \ge 1/2$, (\ref{reclnf4} - \ref{reclnf5}) gives
$$
2 \A_{2,-1}^{-1} (\lambda_2 - \lambda_1) \ge
\frac{C_1 C_2 R(i)}{2} \ge \frac{C_1 C_2}{4} \ge 
\frac{C_1 C_2 \log(L + 1)}{4 L^2},
$$
which proves (\ref{reclnf0}).

So now we assume that for all $i \in \{1,\ldots,[N^{1/4}]\}$ we have $R(i) < S(i)$. This gives $S(i) \ge L(i)/2 \ge 1/2$.

Let us observe that 
\begin{eqnarray}
\nonumber
&&
2 \A_{2,-1}^{-1} (\lambda_2 - \lambda_1) =
\int_{D} \int_{D}
\frac{(f(x) - f(y))^2}{|x - y|^3} \varphi_1(x) \varphi_1(y) \, dx \, dy \ge
\\
\label{reclnf6}
&&
\sum_{i = 1}^{[N^{1/4}]} \sum_{k = 1}^i \sum_{m = p(i,k)}^{q(i,k)-1} \int_{D_{k,m}^i} \int_{D_{k,m + 1}^i}
\frac{(f(x) - f(y))^2}{|x - y|^3} 
\\
&&
\nonumber
\quad \quad \quad \quad \quad \quad \quad \quad\quad \quad \quad \quad \quad
\times \varphi_1(x) \varphi_1(y) \, dx \, dy. 
\end{eqnarray}
Note that 
$$
\sup_{x \in D_{k,m}^i, \, y \in D_{k,m + 1}^i} |x - y|^3 \le ((2 i + 2)^2 + 2^2)^{3/2} \le
(20 i^2)^{3/2} = C_2 i^3.
$$
So by our standard arguments (\ref{reclnf6}) is bounded below by
\begin{eqnarray}
\label{reclnf7}
&&
\sum_{i = 1}^{[N^{1/4}]} \frac{C_1 C_2}{i^3}
\sum_{k = 1}^i \sum_{m = p(i,k)}^{q(i,k) - 1} 
(\mu_{k,m}^i \mu_{k,m + 1}^i)^{-1/2} \int_{D_{k,m}^i} \int_{D_{k,m + 1}^i} 
\\
\label{reclnf7a}
&&
\quad \quad \quad \quad \quad \quad \quad \quad 
\times
(f(x) - f(y))^2 \varphi_1^2(x) \varphi_1^2(y) \, dx \, dy. 
\end{eqnarray}
Note that $|q(i,k) - p(i,k) + 1| \le  2 N/i + 1 \le 3 N/i$. Hence (\ref{reclnf7} - \ref{reclnf7a}) is bounded below by
$$
\sum_{i = 1}^{[N^{1/4}]} \frac{C_1 C_2}{i^3} \left(\frac{i}{3 N}\right)^2 \frac{S(i)}{4},
$$
where $S(i)$ is the expression in (\ref{reclnf3}). We assumed that $S(i) \ge 1/2$. Therefore
$$
2 \A_{2,-1}^{-1} (\lambda_2 - \lambda_1) \ge
\frac{C_1 C_2}{2^3 \cdot 3^2 N^2} \sum_{i = 1}^{[N^{1/4}]} \frac{1}{i} \ge
\frac{C_1 C_2 \log([N^{1/4}] + 1)}{2^3 \cdot 3^2 N^2}.
$$
Note that $([N^{1/4}] + 1)^5 \ge 2 ([N^{1/4}] + 1)^4 \ge 2 N \ge N + 1$. Hence
$\log([N^{1/4}] + 1) \ge (\log(N + 1))/5$. Note also that $L \ge N$ and a function 
$\log(x + 1)/x^2$ is decreasing for $x \ge 1$. Therefore
$$
2 \A_{2,-1}^{-1} (\lambda_2 - \lambda_1) \ge
\frac{C_1 C_2 \log(N + 1)}{2^3 \cdot 3^2 \cdot 5 N^2} \ge
\frac{C_1 C_2 \log(L + 1)}{360 L^2}.
$$
This shows (\ref{reclnf0}) and finishes Case 1.

Now let us consider Case 2. In this case we will show the following lemma.
\begin{lem}
\label{lemcase2}
If $N \ge 16$ and there exist $i \in \{1,\ldots,[N^{1/4}]\}$ such that
\begin{equation}
\label{reclnf8}
\sum_{k = 1}^i \left(\int_{A_k^i} \varphi_1^2 \right)^{-1}
  \left(\int_{A_k^i} f \varphi_1^2 \right)^2 > \frac{1}{2}
\end{equation}
then
$$
2 \A_{2,-1}^{-1} (\lambda_2 - \lambda_1) \ge 2 C_1 C_2 
\left( \frac{1}{256 i^3} - \frac{72}{N}\right).
$$
\end{lem}
Before we come to the proof of this lemma (which is quite technical) let us first show how this lemma implies (\ref{reclnf0}).

We know (Case 2) that $(\ref{reclnf8})$ holds for some $i \in \{1,\ldots,[N^{1/4}]\}$. 
Hence for $N \ge 16$ we have
\begin{equation}
\label{reclnf9}
2 \A_{2,-1}^{-1} (\lambda_2 - \lambda_1) 
\ge \frac{2 C_1 C_2}{256 N} \left( \frac{N}{i^3} - 72 \cdot 256 \right) 
\ge \frac{2 C_1 C_2}{256 N} \left( N^{1/4} - 18432 \right).
\end{equation}
When (say) $N \ge 10^{18}$ then $N^{1/4} \ge 3 \cdot 10^4$ and (\ref{reclnf9}) implies (\ref{reclnf0}).

When $N \le 10^{18}$ we have $\log(L + 1) \le \log(N + 2) \le 42$. Then  Remark \ref{remark1}  implies
$$
2 \A_{2,-1}^{-1} (\lambda_2 - \lambda_1) 
\ge \frac{1}{72 (L + 1)^2} 
\ge \frac{\log(L + 1)}{72 \cdot 4 L^2 \cdot 42}, 
$$
which also gives (\ref{reclnf0}).
\end{proof}

\begin{proof}[Proof of Lemma~\ref{lemcase2}]
Note that if $i = 1$ then the left hand side of (\ref{reclnf8}) equals $0$. So we may and do assume that $i \ge 2$.

In this proof $i \in \{2,\ldots,[N^{1/4}]\}$ is fixed so we will drop $i$ from the notation. We will write $D_{k,m}$ for $D_{k,m}^i$, $A_k$ for $A_k^i$, $p(k)$, $q(k)$ for $p(i,k)$, $q(i,k)$. We will also introduce the following notation
$$
a_{k,m} = \int_{D_{k,m}} f \varphi_1^2 \left(\int_{D_{k,m}}  \varphi_1^2\right)^{-1},
\quad \quad
b_{k,m} = \int_{D_{k,m}}  \varphi_1^2,
$$
$$
a_k = \int_{A_k} f \varphi_1^2 = \sum_{m = p(k)}^{q(k)} a_{k,m} b_{k,m},
\quad \quad
b_k = \int_{A_k}  \varphi_1^2 = \sum_{m = p(k)}^{q(k)}  b_{k,m}.
$$
The condition (\ref{reclnf8}) written in our notation is 
\begin{equation}
\label{case2f0}
\sum_{k = 1}^i \frac{a_k^2}{b_k} > \frac{1}{2}.
\end{equation}

Now we have to estimate $b_k$ from below.
Note that $b_1 + \ldots + b_i = 1$.
Roughly speaking, since $\varphi_1$ is "midconcave", for $N$ large enough $b_1,\ldots,b_i$ have similar values so $b_k \ge c/i$ for $k = 1, \ldots, i$. The following lemma makes the above remark precise.
\begin{lem}
\label{bk}
For $N \ge 16$ and any $k = 1,\ldots,i$ we have
\begin{equation}
\label{bkf1}
b_k \ge 1/(32 i).
\end{equation}
\end{lem}
\begin{proof}
Note that $\sum_{l = 1}^N \int_{D_l} \varphi_1^2 = 1/2$ and $\int_{D_l} \varphi_1^2$ is nonincreasing in $l$ ($l = 1,\ldots,N$) so $\int_{D_1} \varphi_1^2 \ge 1/(2N)$. 

For any $x_2 \in (-1,1)$ the function $x_1 \to \varphi_1(x_1,x_2)$ is concave for $x_1 \in [-L/2,L/2]$ and attains its maximum for $x_1 = 0$. Hence, for any $x_2 \in (-1,1)$ and $x_1 \in [-L/4,L/4]$ we have $\varphi_1(x_1,x_2) \ge \varphi_1(0,x_2)/2$.

Recall that $D_l = ((l-1) L/N, l L/N) \times (-1,1)$. If $l \in [-N/4 + 1,N/4]$ then $(l-1) L/N \ge -L/4$ and $l L/N \le L/4$. It follows that for such $l$ 
$$
\int_{D_l} \varphi_1^2(x_1,x_2) \, dx_1 \, dx_2 \ge
\frac{1}{4} \int_{D_l} \varphi_1^2(0,x_2) \, dx_1 \, dx_2 \ge
\frac{1}{4} \int_{D_1} \varphi_1^2 \ge
\frac{1}{8 N}.
$$ 
Recall that $A_k = \bigcup_{m = p(k)}^{q(k)} D_{k,m} = \bigcup_{m = p(k)}^{q(k)} D_{k + m i}$. 

Let
$C_k = \{m \in \Z: \, -N/4 + 1 \le k + m i \le N/4\}$. For any $m \in C_k$ we have $\int_{D_{k + mi}} \varphi_1^2 \ge 1/(8 N)$ so 
\begin{equation}
\label{bkf2}
b_k = \int_{A_k} \varphi_1^2 = \sum_{m = p(k)}^{q(k)} \int_{D_{k + m i}} \varphi_1^2 \ge \frac{\# C_k}{8 N}, 
\end{equation}
where $\# C_k$ is the number of elements of $C_k$. 
We have
$$
\# C_k \ge \left[\frac{2 [N/4]}{i}\right] \ge \frac{2 ((N/4) - 1)}{i} - 1 =
\frac{N - 4 - 2 i}{2 i}.
$$
We have $i \le N^{1/4}$ and $N \ge 16$ so  it is not difficult to show that $N - 4 - 2 i \ge N/2$. Hence $\# C_k \ge N/(4 i)$. Finally this and (\ref{bkf2}) gives (\ref{bkf1}).
\end{proof}

By (\ref{case2f0}) and Lemma~\ref{bk} we obtain $\sum_{k = 1}^i a_k^2 \ge 1/(64 i)$.

Note that $\sum_{k = 1}^i a_k = \int_{D} f \varphi_1^2 = 0$. Therefore using
Lemma~\ref{partition0} we get
\begin{equation}
\label{case2f1}
\sum_{k = 1}^{i - 1} (a_{k + 1} - a_k)^2 \ge \frac{1}{i^2} \sum_{k = 1}^i a_k^2 \ge \frac{1}{64 i^3}.
\end{equation}
We have
\begin{eqnarray*}
&&
2 \A_{2,-1}^{-1} (\lambda_2 - \lambda_1) =
\int_{D} \int_{D}
\frac{(f(x) - f(y))^2}{|x - y|^3} \varphi_1(x) \varphi_1(y) \, dx \, dy
\\
&&
\ge 2 \sum_{k = 1}^{i - 1} \int_{A_k} \int_{A_{k + 1}}
\frac{(f(x) - f(y))^2}{|x - y|^3} \varphi_1(x) \varphi_1(y) \, dx \, dy. 
\end{eqnarray*}

By our definition of $p(k)$ and $q(k)$ it is easy to notice that $p(k + 1) \le p(k)$ and $q(k + 1) \le q(k)$.

It follows that
\begin{eqnarray*}
&&
\int_{A_k} \int_{A_{k + 1}}
\frac{(f(x) - f(y))^2}{|x - y|^3} \varphi_1(x) \varphi_1(y) \, dx \, dy
\\
&&
\ge
\sum_{m = p(k)}^{q(k+1)} \int_{D_{k,m}} \int_{D_{k + 1,m}}
\frac{(f(x) - f(y))^2}{|x - y|^3} \varphi_1(x) \varphi_1(y) \, dx \, dy.  
\end{eqnarray*}
By our standard arguments this is bounded below by
\begin{eqnarray*}
&&
\sum_{m = p(k)}^{q(k+1)} \frac{C_1 C_2}{b_{k,m}^{1/2} b_{k + 1,m}^{1/2}}
\int_{D_{k,m}} \int_{D_{k + 1,m}}
(f(x) - f(y))^2 \varphi_1^2(x) \varphi_1^2(y) \, dx \, dy 
\\
&&
\ge
C_1 C_2 \sum_{m = p(k)}^{q(k+1)} (a_{k,m} - a_{k + 1,m})^2 b_{k,m}^{1/2} b_{k + 1,m}^{1/2}. 
\end{eqnarray*}
The last inequality follows from the argument which has been already used in the last 3 lines in the proof of Lemma~\ref{partition}.
By Schwarz inequality it is bounded below by 
\begin{equation}
\label{case2f2}
C_1 C_2 \left( \sum_{m = p(k)}^{q(k+1)} (a_{k,m} - a_{k + 1,m}) b_{k,m}^{1/2} b_{k + 1,m}^{1/2} \right)^2
\left( \sum_{m = p(k)}^{q(k+1)}  b_{k,m}^{1/2} b_{k + 1,m}^{1/2} \right)^{-1}.
\end{equation}
We have 
$$
\sum_{m = p(k)}^{q(k+1)}  b_{k,m}^{1/2} b_{k + 1,m}^{1/2}
\le
\left( \sum_{m = p(k)}^{q(k+1)}  b_{k,m} \right)^{1/2}
\left( \sum_{m = p(k)}^{q(k+1)}  b_{k + 1,m} \right)^{1/2}
\le 1.
$$
So (\ref{case2f2}) is bounded below by 
$$
C_1 C_2 \left( \sum_{m = p(k)}^{q(k+1)} (a_{k,m} - a_{k + 1,m}) b_{k,m}^{1/2} b_{k + 1,m}^{1/2} \right)^2.
$$
Now let us denote
$$
R_k = \sum_{m = p(k)}^{q(k+1)} (a_{k,m} - a_{k + 1,m}) b_{k,m}^{1/2} b_{k + 1,m}^{1/2} .
$$
We have $R_k = S_k + T_k + U_k + V_k$, where
$$
S_k = \sum_{m = p(k)}^{q(k)} a_{k,m} b_{k,m} - \sum_{m = p(k + 1)}^{q(k + 1)} a_{k + 1,m} b_{k + 1,m} = a_k - a_{k + 1},
$$
$$
T_k = \sum_{m = p(k)}^{q(k+1)} a_{k,m} b_{k,m}^{1/2} (b_{k + 1,m}^{1/2} - b_{k,m}^{1/2}),
$$
$$
U_k = \sum_{m = p(k)}^{q(k+1)} - a_{k + 1,m} b_{k + 1,m}^{1/2} 
(b_{k,m}^{1/2} - b_{k + 1,m}^{1/2}),
$$
$$
V_k = - \delta_{k,r(i)} a_{k,q(k)} b_{k,q(k)} + 
\delta_{k,i - r(i)} a_{k + 1,p(k + 1)} b_{k + 1,p(k + 1)},
$$
where $\delta_{x,y} = 1$ when $x = y$ and $\delta_{x,y} = 0$ when $x \ne y$. In other words $V_k = 0$ when $k \ne r(i)$ and $k \ne i - r(i)$. In order to see why an extra term $V_k$ appears let us recall the definition of $p(k)$ and $q(k)$. 
We have $q(k) = [N/i]$ for $k = 1, \ldots, r(i)$, $q(k) = [N/i] -
1$ for $k = r(i) + 1,\ldots,i$ and $p(k) = -[N/i]$ for $k = 1,
\ldots,i - r(i)$, $p(k) = -[N/i] -1$ for $k =i - r(i) + 1,\ldots,i$. 
A nontrivial term $V_k$ appears only if $q(k) \ne q(k + 1)$ ($k = r(i)$) or $p(k) \ne p(k + 1)$ ($k = i - r(i)$).

We know that 
$$
(a + b + c + d)^2 \ge (a^2/4) - b^2 - c^2 - d^2, \quad a,b,c,d \in R,
$$
so 
$$
R_k^2 = (S_k + T_k + U_k + V_k)^2 \ge (S_k^2/4) - T_k^2 - U_k^2 - V_k^2.
$$
By (\ref{case2f1}) we obtain 
$$
\sum_{k = 1}^{i - 1} S_k^2 = \sum_{k = 1}^{i - 1} (a_{k + 1} - a_k)^2  \ge \frac{1}{64 i^3}.
$$
We have already obtained that 
\begin{eqnarray}
\nonumber
2 \A_{2,-1}^{-1} (\lambda_2 - \lambda_1) &\ge& 
2 C_1 C_2 \sum_{k = 1}^{i - 1} R_k^2
\\
\nonumber
&\ge& 2 C_1 C_2 \left(\frac{1}{4} \sum_{k = 1}^{i - 1} S_k^2 - \sum_{k = 1}^{i - 1} (T_k^2 + U_k^2 + V_k^2) \right) 
\\
\label{case2f3}
&\ge& 2 C_1 C_2 \left(\frac{1}{256 i^3} - \sum_{k = 1}^{i - 1} (T_k^2 + U_k^2 + V_k^2) \right).
\end{eqnarray}
Now we have to estimate $\sum_{k = 1}^{i - 1} (T_k^2 + U_k^2 + V_k^2)$.

By Schwarz inequality we obtain
\begin{equation}
\label{case2f4}
T_k^2 \le \sum_{m = p(k)}^{q(k+1)} a_{k,m}^2 b_{k,m} |b_{k + 1,m}^{1/2} - b_{k,m}^{1/2}|
\sum_{m = p(k)}^{q(k+1)} |b_{k + 1,m}^{1/2} - b_{k,m}^{1/2}|.
\end{equation}
We have $b_{k,m} = \int_{D_{k,m}} \varphi_1^2 = \int_{D_{k + m i}} \varphi_1^2$. The sequence $\{\left(\int_{D_l} \varphi_1^2\right)^{1/2}\}_{l = -N + 1}^{l = N}$ is unimodal and its maximum is equal to $\left(\int_{D_1} \varphi_1^2\right)^{1/2}$. 

Now there is a very important observation in the proof of this lemma. By the unimodality of this sequence we have
$$
\sum_{m = p(k)}^{q(k+1)} |b_{k + 1,m}^{1/2} - b_{k,m}^{1/2}| \le 2 \left(\int_{D_1} \varphi_1^2\right)^{1/2}.
$$
We also have $|b_{k + 1,m}^{1/2} - b_{k,m}^{1/2}| \le \left(\int_{D_1} \varphi_1^2\right)^{1/2}$. 

On the other hand by Lemma~\ref{phi1est} we know that $||\varphi_1^2||_{\infty} \le 9/L$ and the area $|D_1| = 2 L/N$. Hence $\int_{D_1} \varphi_1^2 \le 18/N$. 

By (\ref{case2f4}) we obtain
$$
T_k^2 \le \frac{36}{N} \sum_{m = p(k)}^{q(k+1)} a_{k,m}^2 b_{k,m}.
$$
Similarly we get
$$
U_k^2 \le \frac{36}{N} \sum_{m = p(k)}^{q(k+1)} a_{k + 1,m}^2 b_{k + 1,m}.
$$

Now we estimate $V_k^2$. Recall that $b_{k,m} \le \int_{D_1} \varphi_1^2 \le 18/N$. We have
\begin{eqnarray*}
&&
V_k^2 \le 2( \delta_{k,r(i)} a_{k,q(k)}^2 b_{k,q(k)}^2 + 
\delta_{k,i - r(i)} a_{k + 1,p(k + 1)}^2 b_{k + 1,p(k + 1)}^2)
\\
&&
\le (36/N) ( \delta_{k,r(i)} a_{k,q(k)}^2 b_{k,q(k)} + 
\delta_{k,i - r(i)} a_{k + 1,p(k + 1)}^2 b_{k + 1,p(k + 1)})
\end{eqnarray*}

It follows that $\sum_{k = 1}^{i - 1} (T_k^2 + U_k^2 + V_k^2)$ is bounded above by
$$
\frac{36}{N} \sum_{k = 1}^{i - 1} \sum_{m = p(k)}^{q(k)} a_{k,m}^2 b_{k,m} + 
\frac{36}{N} \sum_{k = 1}^{i - 1} \sum_{m = p(k + 1)}^{q(k+1)} a_{k + 1,m}^2 b_{k + 1,m}.
$$
Note also that $a_{k,m}^2 b_{k,m} \le \int_{D_{k,m}} f^2 \varphi_1^2$. Hence
$$
\sum_{k = 1}^{i - 1} (T_k^2 + U_k^2 + V_k^2) \le 2 \frac{36}{N} \int_{D} f^2 \varphi_1^2 = 
\frac{72}{N}.
$$
This and (\ref{case2f3}) gives the assertion of the lemma.
\end{proof}

\section{Spectral gap for convex double symmetric domains} \label{sec:convex}

\begin{prop}\label{gapconvexl2}
Let $D\subset [-L,L]\times[-1,1]$ be open, convex and symmetric with respect
to both axis. Assume $(L,0)\in \overline{D}$, $(0,1/2)\in\overline{D}$
and $L=l^2+4$ for some natural number $l\geq 3$.
Then
$$
 \int_D f^2(x) \varphi_1^2(x)\,dx \leq
%tk dodalem 2
2 \cdot 10^9 c_H L^2
 \int_D \!\int_D \frac{(f(x)-f(y))^2}{|x-y|^{2+\alpha}}\,\varphi_1(x)\varphi_1(y)
 \,dx\,dy
$$
for all $f\in L^2(D)$ such that $\int_D f(x)\varphi_1^2(x)\,dx=0$, where
$c_H$ denotes the constant from Corollary~\ref{corHarnack}.
\end{prop}
\begin{proof}
We denote by $D(a,b)$ the set 
$D\cap ((a,b)\times\R)$, or
$D\cap ((a,b]\times\R)$, or
$D\cap ([a,b)\times\R)$. The latter three sets differ
only by a set of a measure zero, thus the ambiguity of
the definition of $D(a,b)$ will be irrelevant.
We put 
$w(D(a,b)) = b-a$, which is the ``width'' of $D(a,b)$,
and 
\begin{eqnarray*}
h(D(a,b))&=&2\inf\{t : (x,t)\in D(a,b) \textrm{ for some $x\in (a,b)$}\},\\
H(D(a,b))&=&2\sup\{t : (x,t)\in D(a,b) \textrm{ for some $x\in (a,b)$}\},
\end{eqnarray*}
the ``heights'' of the set $D(a,b)$.

Let $\mu = \varphi_1^2 \,dx$. We fix an arbitrary $f\in L^2(D,\mu)$ such that
$\int_D f\,d\mu =0$ and
put $F(x,y)=\frac{(f(x)-f(y))^2}{|x-y|^{2+\alpha}}\,\varphi_1(x)
\varphi_1(y)$.

Step 1.
We consider a partition of $D$ into a union of five disjoint sets
$D_1=D(-L,-L+4)$, $D_2=D(-L+4,-L+8)$, $D_3=D(-L+8,L-8)$, 
$D_4=D(L-8,L-4)$ and $D_5=D(L-4,L)$.
Note that by unimodality and symmetry of $\varphi_1^2$ and $w(D_k)$,
the sequence $\mu_k=\int_{D_k} \varphi_1^2\,d\mu$ is also
unimodal.
Thus by Lemma~\ref{partition} we have
\begin{eqnarray}
%tk poniewaz sie nie miescilo pozmienialem troche uklad
&& 2\int_D f^2 \,d\mu \leq
   2\sum_{k=1}^5 \frac{1}{\mu_k}
     \int_{D_k}\!\int_{D_k} (f(x)-f(y))^2\mu(dx)\mu(dy) \nonumber\\
   &+& 
   100 \sum_{k=1}^4 \frac{1}{\mu_k \vee \mu_{k+1}}
   \int_{D_k}\!\int_{D_{k+1}} (f(x)-f(y))^2\mu(dx)\mu(dy) \nonumber\\
   &\leq&
   2(20)^{1+\alpha/2}c_3 L^2\sum_{k=1,k\neq 3}^5
   \int_{D_k}\!\int_{D_k} F(x,y)\,dx\,dy \nonumber\\
   &+& 
   100(68)^{1+\alpha/2} 2 c_3 L^2 \left(
     \int_{D_1}\!\int_{D_2} +  \int_{D_4}\!\int_{D_5} \right)
     F(x,y)\,dx\,dy \nonumber\\
   &+&
   \frac{100}{\mu_2+\mu_3+\mu_4}
   \iint_{(D_2\cup D_3\cup D_4)^2} (f(x)-f(y))^2\,\mu(dx)\mu(dy).
      \label{convexlastint}
\end{eqnarray}
In the above inequality we have used the fact that 
$\sup_E\varphi_1^2 \leq c_3L^2 \int_E \varphi_1^2\,dx$
for $E=D_k$, where $k\neq 3$, or $E=D_1\cup D_2$
or $E=D_4\cup D_5$. It turns out that one may take 
$c_3=9c_H$
%tk tu dodalem pare slow
(This follows from Corollary \ref{corHarnack} an argument of a geometric nature is omitted).
Moreover, $\diam(D_k)\leq \sqrt{20}$
for $k\neq 3$ and $\diam(D_k\cup D_{k+1})\leq \sqrt{68}$
for $k=1,4$.

In this step we have ``cut off'' the ends of $D$, in a sense
that it remains to estimate from above the term in 
(\ref{convexlastint})

Step 2.
We now define a sequence $a_k=l^2 - \sum_{j=1}^k (2j-1)$
for $k=1,2,\ldots, l$, and $a_0=l^2$. Note that $a_l=0$.
We consider a partition of $D(-l^2,l^2)=D_2\cup D_3\cup D_4$
into a union of $2l$ pairwise disjoint sets
$D_{-k}'=D(-a_{k-1},-a_k)$ and $D_k'=D(a_k,a_{k-1})$ for
$k=1,2,\ldots,l$.
Let $\mu_k'=\mu(D_k')$. By a similar token as before, the sequence
$(\mu_{-1}',\mu_{-2}',\ldots,\mu_{-l}',\mu_{l}',\mu_{l-1}',\ldots,\mu_1')$
is unimodal. Thus by Lemma~\ref{partition} and the equality $\mu_k'=\mu_{-k}'$
we have
\begin{eqnarray}
%tk tu znowu sie nie miescilo wiec zmienilem nieco uklad
&&
\frac{1}{\mu(D(-l^2,l^2))} \iint_{D(-l^2,l^2)^2} (f(x)-f(y))^2 \,
\mu(dx)\mu(dy) \nonumber\\
 &\leq &
   2\sum_{k=1}^l \frac{1}{\mu_k'}
   ( \int_{D_k'} \!\int_{D_k'} + \int_{D_{-k}'} \!\int_{D_{-k}'}
   ) 
    (f(x)-f(y))^2 \,\mu(dx)\mu(dy) \label{step2int1}\\
   && +
   16l^2\sum_{k=1}^{l-1} \frac{1}{\mu_{k+1}'\vee \mu_k'}
   ( \int_{D_{-k}'} \!\int_{D_{-k-1}'} + \int_{D_{k+1}'} \!\int_{D_k'}
   ) 
    (f(x)-f(y))^2  \label{step2int2}\\
 && \quad \quad \quad \quad \quad \quad \quad \quad \quad \quad \quad 
\quad \quad \quad \quad \quad \quad \quad \quad  
\times \mu(dx)\mu(dy) 
\nonumber \\
   && +
   16l^2 \frac{1}{\mu_l'}
   \int_{D_{-l}'} \!\int_{D_l'}  
    (f(x)-f(y))^2 \,\mu(dx)\mu(dy).\label{step2int3}
\end{eqnarray}

Step 3.
We will now show how to deal with the integral (\ref{step2int2}), i.e.,
\begin{eqnarray}
I_k &=&
\frac{1}{\mu_{k+1}'\vee \mu_k'}
    \int_{D_{k+1}'} \!\int_{D_k'}
       (f(x)-f(y))^2 \,\mu(dx)\mu(dy) \nonumber\\
  &\leq&
\frac{1}{\mu_{k+1}'+ \mu_k'}
    \iint_{(D_{k+1}' \cup D_k')^2}
       (f(x)-f(y))^2 \,\mu(dx)\mu(dy). \label{step3Ik}
\end{eqnarray}
We have $w = w(D_{k+1}'\cup D_k')=a_{k-1}-a_{k+1}=4k$.
Let $h=h(D_{k+1}'\cup D_k')$ and
$N=\left[ \frac{2k+1}{h} \right]$.
We divide $D_{k+1}'\cup D_k'$ into a union of sets
$$
E_j=D(a_{k+1}+(j-1)w/N, a_{k+1}+jw/N),\quad
j=1,2,\ldots,N,
$$
of equal width $4k/N$ and
apply Lemma~\ref{partition} to such $E_j$. We
obtain
\begin{eqnarray*}
I_k &\leq&
 4\sum_{j=1}^N \frac{1}{\mu(E_j)}
  \int_{E_j}\!\int_{E_j} (f(x)-f(y))^2 \,\mu(dx)\mu(dy) \\
  && +
 8N^2\sum_{j=1}^{N-1} \frac{1}{\mu(E_j)\vee \mu(E_{j+1})}
  \int_{E_j}\!\int_{E_{j+1}} (f(x)-f(y))^2 \,\mu(dx)\mu(dy).
\end{eqnarray*}
Note that
$$
 \dist((D_{k+1}'\cup D_k'), (L,0)) = 4+(k-1)^2 \geq 2k = w(D_{k+1}'\cup D_k')/2,
$$
thus by convexity of $D$ we have
$$
  \frac{H(D_{k+1}'\cup D_k')}{h(D_{k+1}'\cup D_k')} \leq
  \frac{ \dist((D_{k+1}'\cup D_k'), (L,0)) +w(D_{k+1}'\cup D_k')}
  {\dist((D_{k+1}'\cup D_k'), (L,0))} \leq 3.
$$
Hence $h\leq H(E_j)\leq 3h$. Moreover,
$$
 w(E_j\cup E_{j+1})=\frac{8k}{N} \leq 8h \frac{k}{2k+1-h} \leq 4h
$$
and
$$
 w(E_j\cup E_{j+1})=\frac{8k}{N} \geq h \frac{8k}{2k+1} \geq \frac{8}{3}h,
$$
This means that if $S=a_{k+1}+(j-1)w/N$, then $B((S,0),h)\subset D$ and
$|B((S,0),h/4)\cap E_j| = |B((S,0),h/4)|/2 = \pi h^2/32$. Thus by Harnack inequality
(Corollary \ref{corHarnack}) we obtain
$$
 \sup_{E_j} \varphi_1^2  \leq \frac{32c_H}{\pi h^2} 
   \int_{E_j} \varphi_1^2 \,dx,
$$
the same bound as above holds also for $E_j\cup E_{j+1}$ in place of $E_j$.
%bd Dodaje taka uwage:
Note that we take $h/4$ above as the radius of the ball because of the assumption
concerning inner radius in Corollary \ref{corHarnack}.

We have
$$
\diam(E_j) \leq \diam(E_j\cup E_{j+1}) \leq
 (w(E_j\cup E_{j+1})^2+H(E_j\cup E_{j+1})^2)^{1/2} \leq 5h.
$$
Hence 
$$ 
 \frac{1}{\mu(E_j)\vee \mu(E_{j+1})}
  \int_{E_j}\!\int_{E_{j+1}} (f(x)-f(y))^2 \,\mu(dx)\mu(dy)
$$$$ \leq
  (5h)^{2+\alpha} \frac{32c_H}{\pi h^2} 
  \int_{E_j}\!\int_{E_{j+1}} F(x,y) \,dx\,dy
$$
and a similar bound holds for the integral over $E_j\times E_j$.
Moreover, $N^2\leq (2k+1)^2/h^2 \leq 9k^2/h^2$. Thus
\begin{eqnarray}
I_k &\leq&
4 (5h)^{2+\alpha} \frac{32c_H}{\pi h^2} 
 \sum_{j=1}^N 
  \int_{E_j}\!\int_{E_j} F(x,y)\,dx\,dy \nonumber\\
  && +
 \frac{72k^2}{h^2}\cdot (5h)^{2+\alpha} \frac{32c_H}{\pi h^2} 
  \sum_{j=1}^{N-1} 
  \int_{E_j}\!\int_{E_{j+1}} F(x,y)\,dx\,dy\nonumber\\
  &\leq&
%tk 2304
 \frac{2304\cdot 5^{2+\alpha} c_H}{\pi}\cdot k^2 h^{\alpha-2} 
  \iint_{(D_{k+1}'\cup D_k')^2} F(x,y)\,dx\,dy. \label{estIk}
\end{eqnarray}
We have
$\dist((D_{k+1}'\cup D_k'), (L,0)) = 4+(k-1)^2 \geq k^2/2$,
thus that by convexity of the set $D$ and the assumptions
 $(0,1/2)\in \overline{D}$, $(L,0)\in \overline{D}$
 we obtain
$$
 \frac{h/2}{k^2/2} \geq \frac{1/2}{L}.
$$
When $\alpha \ge 1$ we get $k^2 h^{\alpha - 2} \le 2 L h h^{\alpha - 2} \le 2 L$. When $\alpha < 1$ we get $k^2 h^{\alpha - 2} \le k^2 (2L/k^2)^{2 - \alpha} = 2^{2 - \alpha} L (L/k^2)^{1 - \alpha} \le 2^{2 - \alpha} L$. For any $\alpha \in (0,2)$ we have $k^2 h^{\alpha - 2} \le \max(2,2^{2 - \alpha})L$.
We combine it with (\ref{estIk}) and finally obtain
\begin{equation}\label{estIk2}
I_k \leq
%tk 2880000
 \frac{2880000 c_H}{\pi}\cdot L 
  \iint_{(D_{k+1}'\cup D_k')^2} F(x,y)\,dx\,dy. 
\end{equation}

We should also estimate from above the integral (\ref{step2int3})
over $D_{-l}\times D_l$. This may be done in a similar way
as the integrals $I_k$ above.
%tk wydaje mi sie ze podzial powinien byc inny
%tk na N=[2l/h] czesci bo h = h(D_{-l}\cup D_l) niekoniecznie jest 1
%tk ale po prostu te uwagi jak podzilic w ogole usunalem
%tk  wydaje mi sie ze nie sa potrzebne
%; we may for example divide
%$D_{-l}\times D_l = D(-(2l-1),2l-1)$ into a union of $2l-1$
%disjoint sets $E_j=D(b_j,b_j')$ of equal width $2$, then
%we proceed analogously. 
We obtain a similar estimate as
(\ref{estIk2}) with slightly smaller constant, we omit the details.

To estimate (\ref{step2int1}) we may see that in (\ref{step3Ik})
we have in fact estimated $I_k$ from above by an integral
over $(D_{k+1}' \cup D_k')^2$. Thus 
a similar estimation
as in (\ref{estIk2}) holds also for the integrals in (\ref{step2int1}).

We finally obtain
%tk 4/L i 2880000
\begin{eqnarray*}
 2\int_D f^2\,d\mu &\leq&
 \left(  18\cdot (20)^{1+\alpha/2} + 1800 \cdot (68)^{1 + \alpha/2} +
100(\frac{4}{L}+32)
\cdot \frac{2880000}{\pi}
  \right) \\
&\times& c_H L^2
 \iint_{D^2} F(x,y)\,dx\,dy
\end{eqnarray*}
and the proposition follows.
\end{proof}

\begin{proof}[Proof of Theorem~\ref{gapconvexset}.]
By scaling of eigenvalues (\ref{gapscaling}) it is sufficient to consider domains $D$ such that $[-L,L] \times [-1,1]$, $L \ge 1$ is the smallest rectangle (with sides parallel to the coordinate axes) containing $D$ and to show that for such domains the following inequality holds:
$$
 2\mathcal{A}_{2,-\alpha}^{-1} (\lambda_2-\lambda_1)
  \geq  \frac{C}{L^2},
$$
where $C$ is the same as in (\ref{constantC}).

At first assume that $L \ge 13$. For any natural number $l\geq 3$ we have
$$
 \frac{(l+1)^2+4}{l^2+4} \leq \frac{20}{13},
$$
thus there exists $\beta\in[13/20,1]$ such that
 $\beta L=l^2+4$ for some natural number $l\geq 3$.
By (\ref{gapscaling}) we have
$$
 2\mathcal{A}_{2,\alpha}^{-1} (\lambda_2(D)-\lambda_1(D))=
 \beta^{\alpha}\cdot 2\mathcal{A}_{2,\alpha}^{-1} (\lambda_2(\beta D)-\lambda_1(\beta D)).
$$
Note that $\beta D$ satisfies assumptions of Proposition~\ref{gapconvexl2} (in particular $(0,1/2) \in \beta D$). Hence by Theorem~\ref{variational} and Proposition~\ref{gapconvexl2} we obtain
\begin{eqnarray*}
 2\mathcal{A}_{2,\alpha}^{-1} (\lambda_2(D)-\lambda_1(D)) 
 &=&
 2\mathcal{A}_{2,\alpha}^{-1} \beta^{\alpha} (\lambda_2(\beta D)-\lambda_1(\beta D)) \\
 &\ge&
 \frac{\beta^{\alpha}}{10^9 c_H (\beta L)^2} \ge 
 \frac{1}{10^9 c_H L^2}.
\end{eqnarray*}

What remains is to consider the case $L \le 13$. 

Note that $B((0,0),1/\sqrt{2}) \subset D$ so by Corollary \ref{corHarnack} $\varphi_1^2$ satisfies Harnack inequality on 
%bd $B = B((0,0),1/(2 \sqrt{2}))$, -- zmiejszam promien z uwagi na zalozenie o Inr(D) w nier. Harnacka
 $B = B((0,0),1/(4 \sqrt{2}))$, 
in particular $\varphi_1^2(0,0) \le c_H \varphi_1^2(x)$, $x \in B$. Of course $\sup_{x \in D} \varphi_1^2(x) = \varphi_1^2(0,0)$. We have
$$
c_H = c_H \int_D \varphi_1^2 \ge c_H \int_B \varphi_1^2 \ge \varphi_1^2(0,0) |B|,
$$
which gives $\sup_{x \in D} \varphi_1^2(x) \le c_H |B|^{-1} = 32 c_H/\pi$.

We also have $\diam(D) < 28$. Let $f = \varphi_2/\varphi_1$. By Theorem~\ref{variational} we have
\begin{eqnarray*}
&& 2 \mathcal{A}_{2,\alpha}^{-1} (\lambda_2 - \lambda_1) 
= \int_D \int_D \frac{(f(x) - f(y))^2}{|x - y|^{2 + \alpha}} \varphi_1(x) \varphi_1(y) \, dx \, dy \\
&& \ge \frac{\pi}{32 c_H 28^{2 + \alpha}} \int_D \int_D (f(x) - f(y))^2  \varphi_1^2(x) \varphi_1^2(y) \, dx \, dy \\
&& = \frac{\pi}{32 c_H 28^{2 + \alpha}} 2 \int_D f^2(x) \varphi_1^2(x) \, dx
= \frac{\pi}{16 c_H 28^{2 + \alpha}}.
\end{eqnarray*} 
\end{proof}

{\bf{Acknowledgments.}} The second name author thanks very much Rodrigo Ba\~nuelos for many interesting discussions on the subject of the paper.


\begin{thebibliography}{999}

\bibitem{B} R. Ba\~nuelos,
\emph{Intrinsic ultracontractivity and eigenfunction estimates for
  Schr{\"o}dinger operators}, J. Funct. Anal. 100 (1991), 181-206.

\bibitem {BK1}
R. Ba\~nuelos, T. Kulczycki,
\emph{The Cauchy process and the Steklov problem},
J. Funct. Anal. 211(2) (2004), 355-423.

\bibitem{BK2}
R. Ba\~nuelos, T. Kulczycki,
\emph{Eigenvalue gaps for the Cauchy process and a Poincare inequality},
J. Funct. Anal. 234(1) (2006) 199-225.

\bibitem{BK3}
R. Ba\~nuelos, T. Kulczycki,
\emph{Spectral gap for the Cauchy process on convex, symmetric domains},
Comm. Partial Differential Equations 31(12)(2006), to appear.


\bibitem{BKM}
R{.} Ba\~nuelos, T{.} Kulczycki, P{.} J{.} M\'endez-Hern\'andez,
\emph{On the shape of the ground state eigenfunction for stable
processes},  Potential Anal. 24(3) (2006), 205-221. 

\bibitem{BLM}
R{.} Ba\~nuelos,
R{.} Lata\l{}a, P{.} J{.} M\'endez-Hern\'andez,
\emph{A
Brascamp-Lieb-Luttinger-type
inequality and applications to symmetric
stable processes},
Proc{.} Amer{.} Math{.} Soc{.} 129(10) (2001),
2997--3008 (electronic).

\bibitem{BM}
R{.} Ba\~nuelos, P{.} J{.} M\'endez-Hern\'andez,
\emph{Sharp inequalities for heat kernels of
Schr{\"o}dinger operators and
applications to spectral gaps},
J{.}
Funct{.} Anal{.} 176(2) (2000), 368--399.

\bibitem{BG2}
R{.}M{.} Blumenthal and R{.}K{.}
Getoor,
\emph{The asymptotic distribution of the eigenvalues for a
class of Markov
operators}
Pacific J. Math. 9 (1959),
399--408.

\bibitem{BGR}
R.M. Blumenthal, R.K. Getoor, D.B. Ray \emph{On the distribution of
  first hits for the symmetric stable processes},
Trans. Amer. Math. Soc. 99 (1961) 540-554.

\bibitem{BB} K. Bogdan and T. Byczkowski \emph{Potential theory for
    the $\alpha$-stable Schr\"odinger operator on bounded Lipschitz
    domains}, Studia Math. 133(1) (1999), 53-92.

\bibitem{BK} K. Burdzy, T. Kulczycki \emph{Stable processes have
    thorns}, Ann. Probab. 31(1) (2003), 170-194.

\bibitem{CS1}
Z{.}Q{.} Chen and R{.}
Song,
\emph{Intrinsic ultracontractivity and conditional gauge for
symmetric  stable processes},
J{.} Funct{.} Anal{.} 150(1) (1997),
204--239.


\bibitem{CS2}
Z{.}Q{.} Chen and R{.}
Song,
\emph{Intrinsic ultracontractivity, conditional lifetimes and
   conditional gauge for symmetric stable processes on rough  domains},
Illinois J{.} Math{.} 44(1) (2000),
138--160.

\bibitem{CS3}
Z{.}Q{.} Chen and R{.}
Song,
\emph{Two sided eigenvalue estimates for subordinate Brownian motion in bounded domains},
J. Funct. Anal. 226 (2005), 90-113.

\bibitem{CS4}
Z{.}Q{.} Chen and R{.}
Song,
\emph{Continuity of eigenvalues for subordinate processes in domains},
Math. Z. (2005), (to appear).


\bibitem{Da1} E{.}B{.} Davies, Heat Kernels and
Spectral Theory, Cambridge
University Press,
Cambridge,
1989.

\bibitem{Da} B{.} Davis,
\emph{On the
spectral gap for fixed membranes},
Ark{.} Mat{.} 39(1) (2001),
65--74.

\bibitem{Db} R. D. DeBlassie,
\emph{Higher order PDEs and symmetric stable processes}
Probab. Theory Related Fields 129 (2004), 495--536.

\bibitem{DbM} R. D. DeBlassie and P{.} J{.} M\'endez-Hern\'andez,
\emph{$\alpha$--continuity properties of symmetric $\alpha$--stable process},  Preprint. 

%tk nowe cytowanie
\bibitem{G} R{.} K{.}
Getoor,
\emph{First passage times for symmetric stable processes in space},
Trans. Amer. Math. Soc. 101 (1961)
75-90.

\bibitem{G1} R{.} K{.}
Getoor,
\emph{Markov operators and their associated
semi-groups},
Pacific J{.} Math{.} 9 (1959)
449--472.

\bibitem{IW} N. Ikeda, S. Watanabe \emph{On some relations between the
  harmonic measure and the Levy measure for a certain class of Markov 
processes}, J. Math. Kyoto Univ. 2 (1962), 79-95.

\bibitem{K}
T{.} Kulczycki,
\emph{Intrinsic
ultracontractivity for symmetric stable processes},
Bull{.} Polish
Acad{.} Sci{.} Math{.} 46(3) (1998),
325--334.

\bibitem{L}
J{.} Ling, \emph{A lower bound for the gap between the first two eigenvalues of
Schr{\"o}dinger operators on convex domains in $S\sp n$ or $R\sp n$},
Michigan Math{.} J{.} 40(2) (1993), 259--270.

\bibitem{M}
P{.} J{.} M\'endez-Hern\'andez, \emph{Brascamp-Lieb-Luttinger Inequalities for Convex Domains of Finite Inradius},
Duke Math. Journal 113 (2002), 93--131. 

\bibitem{SWYY}
I{.} M{.} Singer, B{.} Wong,
S{.}-T{.}
Yau, S{.} S{.}-T{.} Yau,
\emph{An estimate of the gap of
the first
two eigenvalues in
the Schr{\"o}dinger operator},
Ann{.}
Scuola
Norm{.} Sup{.} Pisa Cl{.} Sci{.} (4) 12(2) (1985),  319--333.

\bibitem{Sm} R. G. Smits,  
\emph{Spectral gaps and rates to equilibrium for diffusions in convex
  domains}  Michigan Math. J. 43(1) (1996), 141-157. 

\bibitem{YZ} Q. Yu and J. Q. Zhong,
\emph{Lower bounds of the gap between the first and second eigenvalues
of the Schr{\"o}dinger operator}, Trans. Amer. Math. Soc. 294 (1986), 341-349.

\bibitem{Z} V. M. Zolotarev, \emph{Integral transformations of
 distributions and estimates of parameters of multidimensional
 spherically symmetric stable laws},
 in: Contributions to Probability, Academic Press, New York, (1981), 283-305.

\end{thebibliography}
\end{document}